\documentclass[a4paper]{amsart}
\usepackage{color}
\usepackage[latin1]{inputenc}
\usepackage[T1]{fontenc}
\usepackage{amsfonts}
\usepackage{amssymb}
\usepackage{amsmath}
\usepackage{amsthm}
\usepackage{stmaryrd}
\usepackage{enumerate}
\usepackage{multirow}
\usepackage{graphicx}
\usepackage{setspace}
\usepackage{graphicx,type1cm,eso-pic,color}
\usepackage{pstricks}
\usepackage{float}
\newtheorem{theorem}{Theorem}[section]
\newtheorem{lemma}[theorem]{Lemma}
\newtheorem{proposition}[theorem]{Proposition}

\newtheorem{assumption}[theorem]{Assumption}
\newtheorem{remark}[theorem]{Remark}
\newtheorem{example}[theorem]{Example}
\begin{document}
\setlength\arraycolsep{2pt}
\title[Fixed-Point Estimation of the Drift in SDE Driven by Rough Fractional Noise]{Fixed-Point Estimation of the Drift Parameter in Stochastic Differential Equations Driven by Rough Multiplicative Fractional Noise}
\author{Chiara AMORINO$^{\star}$}
\email{chiara.amorino@upf.edu}
\author{Laure COUTIN$^{\dag}$}
\email{laure.coutin@math.univ-toulouse.fr}
\author{Nicolas MARIE$^{\diamond}$}
\email{nmarie@parisnanterre.fr}
\date{}
\maketitle
\noindent
$^{\star}$Universitat Pompeu Fabra, Barcelona, Spain.\\
$^{\dag}$Institut de Math\'ematiques de Toulouse, Universit\'e Paul Sabatier, Toulouse, France.\\
$^{\diamond}$Laboratoire Modal'X, Universit\'e Paris Nanterre, Nanterre, France.
%


%
\begin{abstract}
We investigate the problem of estimating the drift parameter from $N$ independent copies of the solution of a stochastic differential equation driven by a multiplicative fractional Brownian noise with Hurst parameter $H\in (1/3,1)$. Building on a least-squares-type object involving the Skorokhod integral, a key challenge consists in approximating this unobservable quantity with a computable fixed-point estimator, which requires addressing the correction induced by replacing the Skorokhod integral with its pathwise counterpart. To this end, a crucial technical contribution of this work is the reformulation of the Malliavin derivative of the process in a way that does not depend explicitly on the driving noise, enabling control of the approximation error in the multiplicative setting. For the case $H\in (1/3,1/2]$, we further exploit results on two-dimensional Young integrals to manage the more intricate correction term that appears. As a result, we establish the well-posedness of a fixed-point estimator for any $H\in (1/3,1)$, together with both an asymptotic confidence interval and a non-asymptotic risk bound. Finally, a numerical study illustrates the good practical performance of the proposed estimator.
\end{abstract}
%


%
\section{Introduction}\label{section_introduction}
In this work, we consider the stochastic differential equation (SDE)
\begin{equation}\label{main_equation}
X_t = x_0 +
\theta_0\int_{0}^{t}b(X_s)ds +\int_{0}^{t}\sigma(X_s)dB_s,
\quad t\in [0,T],
\end{equation}
where $T > 0$ is fixed, $x_0\in\mathbb R$, and $B$ denotes a fractional Brownian motion (fBm) with Hurst parameter $H\in (1/3,1)$. The map $(b,\sigma)$ belongs to $C^{[1/H] + 1}(\mathbb R;\mathbb R^2)$ and has bounded derivatives. The integral with respect to $B$ is understood pathwise, in the sense of Young for $H\in (1/2,1)$, and in the sense of rough paths for $H\in (1/3,1/2]$. Under these conditions, Equation (\ref{main_equation}) admits a unique solution $X$.

Our objective is to estimate the drift parameter $\theta_0$ based on continuous observations of $N$ i.i.d. copies of the equation above, in the asymptotic regime where $T$ is fixed and the number of observed SDEs $N$ tends to infinity.

Stochastic differential equations are widely employed to model real-world phenomena thanks to their versatility. Classical applications are found in finance \cite{LL07}, biology \cite{RICCIARDI77}, neurosciences \cite{HOLDEN76}, and economics \cite{BERGSTROM90}, as well as in physics \cite{PAPANICOLAOU95} and mechanics \cite{KUSHNER67}. In pharmacology, SDEs are commonly used to model variability in biomedical experiments, as reviewed in \cite{AP97, DS13, RAMANATHAN99a, RAMANATHAN99b, WLFWW13}.

Traditionally, inference procedures were designed for a single SDE, either over a fixed time interval or in the limit as $T\rightarrow\infty$. However, with the advent of big data, statisticians increasingly have access to large amounts of data generated by different individuals. This motivates an asymptotic regime where one observes many independent copies over a finite time horizon, each copy corresponding to a different individual. Such a framework has found applications in biology and the social sciences \cite{ALBI19, DEGOND18, MT14, NPT10}, in data-science (notably for neural networks) \cite{CB18, BDFS20, RVE22, SS20}, in optimization \cite{CJLZ21, GP21, TOTZECK21}, and in MCMC methods \cite{DELMORAL13}. For two detailed motivating examples of i.i.d. SDEs driven by standard Brownian motion, we refer to Section 2 of \cite{CM25}.

In recent years, statistical inference from copies of diffusion processes driven by classical Brownian motion - and in particular, estimation of the drift function - has been extensively investigated. See for instance \cite{CG20, MR23, DDM21} for the nonparametric setting, and \cite{DDM20} for the parametric one. The literature is even broader when dependence among copies is allowed, giving rise to observations from diffusions with correlated driving signals (see \cite{CM23}) or interacting particle systems. In this last setting, drift parameter estimation has been thoroughly studied, see \cite{DMH22, BPP23, BPZ24, ABPP23Z24} for nonparametric approaches and \cite{AHPP23, SKPP21} for parametric ones.

In the case of independent copies, extensions to privacy-constrained estimation procedures have been considered in \cite{AGH25}, while replacing Brownian motion with L\'evy processes is studied in \cite{HM24}. All of this confirms that statistical inference from i.i.d. copies of SDEs is a very active area of research of broad interest to the statistics community.

In the classical SDE framework, replacing standard Brownian motion with fractional Brownian motion has proven highly effective for modeling certain real-world phenomena. Indeed, fractional Brownian motion has emerged as a powerful tool for capturing processes driven by long-range dependence and self-similar noise. It has found applications across diverse fields, including finance \cite{CV12, CCR12, FLEMING11, GH04, ROSTEK09}, geophysics \cite{EFTAXIAS08, HURST51}, traffic modeling \cite{CAGLAR04, NORROS95, WTLW95}, and medical research \cite{LDH00}. These developments have naturally stimulated interest in statistical inference for models driven by fBm.

For example, systems of $N$ (dependent) SDEs have recently been employed to model turbulent kinetic energy \cite{BJMR22}. However, empirical evidence suggests that non-Markovian driving noise may offer a more accurate description of such phenomena \cite{FOBWL15, LSEO17}, thus motivating the replacement of standard Brownian motion with fractional Brownian motion. Moreover, \cite{HJMV25} studies the growth of $N$ axon cells in vertebrate brains using the fractional modeling approach, highlighting the relevance and interest of this topic within the scientific community.

It is important to note that statistical inference procedures must be adapted when replacing standard Brownian motion with fractional Brownian motion. In the classical Brownian case, estimation techniques rely crucially on the Markov and semimartingale properties of the process - properties that fractional Brownian motion lacks.

In the context of drift parameter estimation for SDEs driven by fractional Brownian motion, the earliest methods were developed for one long-term observation of the stationary ergodic solution to Equation (\ref{main_equation}), whose existence and uniqueness are guaranteed by a dissipativity condition on the drift function. In this framework, the most common estimation approaches in the literature are the maximum likelihood estimator (MLE) and the least squares estimator (LSE). Notably, for Brownian motion, these two estimators coincide, while for fractional Brownian motion this equivalence no longer holds (see \cite{HN10, HNZ19, KL01, TV07}).

The literature is broadly divided between these two estimation strategies. For instance, MLE-based approaches can be found in \cite{KMR17}, \cite{KL01} and \cite{TV07}, while the LSE has been used in works such as \cite{HN10}, \cite{NT14} and \cite{HNZ19}.

More recently, estimators computed from multiple copies of $X$ driven by fBm, observed over short time intervals, have also been investigated. For example, \cite{DDLW21} analyzes a maximum likelihood estimator, although without theoretical guarantees, while \cite{MARIE25a} develops a procedure based on the LSE, proposing a fixed-point approach to estimate $\theta_0$ when $H\in (1/2,1)$ and with constant $\sigma$. Then, \cite{ANS25} further extends this methodology to interacting particle systems driven by additive fractional noise.

Let us say a few words about the relationship between the two aforementioned observation schemes. As a premise, it is worth acknowledging that the framework based on $N$ independent observations was pioneered for rough differential systems in \cite{PL11}. While this multi-trajectory approach is increasingly natural in contexts where abundant data is available, its relationship with the long-time observation of a single trajectory is non-trivial.

Since standard Brownian motion is a Gaussian process with independent and stationary increments, from one long-time observation of a recurrent diffusion process $Z$, one can always construct $N$ independent copies of $Z_{|[0,T]}$ (see e.g. Marie \cite{MARIE25b}, p. 3-4). This is no longer true from one long-time observation of a fractional diffusion $Y$ because the increments of the fBm $B$ are not independent when $H\neq 1/2$. However, at least when $Y$ is a non-autonomous linear fractional diffusion, since the correlation between $B_{t +\Delta} - B_{\Delta}$ and $B_t$ ($t\in\mathbb R_+$) decreases down to zero when $\Delta\rightarrow\infty$, one can construct $N$ copies $Y^1,\dots,Y^N$ of $Y_{|[0,T]}$ from $(Y_s)_{s\in\mathbb R_+}$ such that
\begin{displaymath}
\sup_{s\in [0,T]}
|\mathbb E(Y_{s}^{i}Y_{s}^{k})|\leqslant\mu_{i,k}(\Delta)
\quad {\rm and}\quad
\lim_{\Delta\rightarrow\infty}\mu_{i,k}(\Delta) = 0,
\quad i,k\in\{1,\dots,N\},
\end{displaymath}
and then establish risk bounds - depending on $\Delta$ - on some estimators of the drift function involved in the definition of $Y$ (see Marie \cite{MARIE26}, Section 3.2 and Corollary 3). This is out of the scope of our paper, but such a strategy could be investigated in order to construct copies - with controlled correlations - of $\overline X_{|[0,T]}$ from one long-time observation of the stationary solution $\overline X$ of Equation (\ref{main_equation}) under a dissipativity condition on $b$.

In this work, we follow this second line of research, choosing to base our estimation procedure on an unobservable quantity, $\widehat\theta_N$, inspired by a least-squares-type estimator. However, this quantity cannot be considered as a genuine estimator, because it depends on the data through a Skorokhod integral, rendering it unobservable. Consequently, borrowing the playful terminology introduced in \cite{ANS25}, we will refer to this quantity as a {\it fakestimator} in the sequel. It can be expressed as
\begin{displaymath}
\widehat\theta_N =
\left(\sum_{i = 1}^{N}\int_{0}^{T}b(X_{s}^{i})^2ds\right)^{-1}
\left(\sum_{i = 1}^{N}\int_{0}^{T}b(X_{s}^{i})\delta X_{s}^{i}\right),
\end{displaymath}
where $N\in\mathbb N^*$, $X^i =\mathcal I(x_0,B^i)$ for every $i\in\{1,\dots,N\}$, $B^1,\dots,B^N$ are independent copies of $B$, $\mathcal I$ denotes the solution map for Equation (\ref{main_equation}), and the stochastic integral is taken in the sense of Skorokhod.

Moving from this fakestimator to a computable estimator involves replacing the Skorokhod integral with a corrected pathwise one. In the case $H > 1/2$, one can exploit the classical relationship between Skorokhod and pathwise integrals, where the Malliavin derivative appears in the correction term. However, since we consider the multiplicative setting, the Malliavin derivative depends on the fBm explicitly in a complicated way, making it very challenging to manage, and consequently to control the error arising when transitioning from the fakestimator to a fixed-point estimator. To overcome this difficulty, our strategy relies on a refined argument that allows us to reformulate the Malliavin derivative of the process so that it can be expressed in a form that does not explicitly depend on $B$ (see Proposition \ref{Malliavin_derivative_X} for details).

Moreover, we also study the case $H\in (1/3,1/2]$, for which this change of integral is even more challenging. In this regime, one can first employ an elegant result from \cite{ST22} to move from the Skorokhod integral to the pathwise one, obtaining a correction term involving the Malliavin derivative trough a 2D Young integral. However, this correction term is difficult to evaluate. Crucially, to obtain the main results for $H\in (1/3,1/2]$, we rely on our Proposition \ref{2D_Young_integral_R}, which allows us to handle the rectangular integrals and recover an expression similar to the case $H > 1/2$. We believe this result is not only essential for our analysis, but is also of independent interest as it could be applied in other contexts of statistical inference for fBm-driven processes with $H\in (1/3,1/2]$.

Thanks to this result, we can move from the Skorokhod integral to the pathwise integral with a computable correction term depending on $\theta_0$, the parameter we wish to estimate. This justifies introducing a fixed-point estimator, denoted by $\overline\theta_N$.

Let us remark that the case $H\in (1/3,1/2]$ is particularly challenging, as confirmed by the scarcity of literature in this context. The same applies to statistical inference with multiplicative fractional noise. Indeed, among all the previously cited works on statistical inference for fractional-noise-driven systems, only \cite{HNZ19} addresses the multiplicative case, in the asymptotic regime with $N = 1$ and $T\rightarrow\infty$.

Note that, by the change-of-variable formula for the pathwise integral, applying a Lamperti-type transform $Y =\Sigma\circ X$, where $\Sigma$ is a primitive function of $1/\sigma$, yields the following stochastic differential equation with additive noise:
\begin{displaymath}
Y_t = Y_0 +
\int_{0}^{t}\left(\frac{b\circ\Sigma^{-1}}{\sigma\circ\Sigma^{-1}}\right)(Y_s)ds + B_t,
\quad t\in [0,T].
\end{displaymath}
This naturally suggests two distinct estimation strategies. The first is a transform-then-apply approach: for $H\in (1/2,1)$, one could transform the original data into $\Sigma\circ X^1,\dots,\Sigma\circ X^N$ and apply the known results for additive fractional noise from \cite{MARIE25a}. However, although $\sigma$ is assumed to be known, there is generally no closed-form expression for the primitive $\Sigma$ or its inverse $\Sigma^{-1}$. Consequently, this method requires numerical integration at every observation point, which is computationally expensive and introduces approximation errors that inevitably propagate into the final estimation. This renders the transform approach highly unsatisfactory from a statistical perspective.

To bypass these limitations, we opt for a direct method based on the model with multiplicative noise. We establish our theoretical results directly on the estimator $\overline\theta_N$ computed from the observations of $X_{|[0,T]}$. The primary gain of this approach is that it completely circumvents the evaluation of $\Sigma$, yielding a robust, purely data-driven procedure free from intermediate numerical integration errors. Furthermore, it establishes a unified methodological framework capable of handling the rougher regime $H\in (1/3,1/2]$. The technical cost of this direct approach, however, is substantial: working directly on the model with multiplicative noise precludes the use of standard additive formulations and demands additional Malliavin calculus-based investigations.

Moreover, note that for an additive fractional noise with Hurst parameter $H \in (1/3,1/2]$, the fixed-point estimation strategy was already introduced in \cite{MARIE25a}, but without any theoretical guarantees. In contrast, thanks to the probabilistic tools developed in this work, we are able to show that, for any $H \in (1/3, 1]$, the fixed-point estimator is well-defined, and we can establish both an asymptotic confidence interval and a non-asymptotic risk bound for $\overline\theta_N$. We emphasize that this constitutes a major step forward in the statistical inference for processes driven by fractional Brownian motion, as it simultaneously addresses (and solves!) two major challenges: the multiplicative noise setting and the regime $H < 1/2$.

The results provided in this paper not only offer valuable insights for tackling related challenges - for instance, nonparametric estimation of the drift or parametric inference under high-frequency observations in the presence of multiplicative noise and/or with $H < 1/2$ - but also validate their practical relevance. Indeed, we conclude our analysis with a numerical study demonstrating that our theoretical fixed-point estimator also performs well in practice.\\

The outline of the paper is as follows. Section \ref{section_probabilistic_preliminaries} introduces the probabilistic tools mentioned above. All theoretical guarantees for our fixed-point estimator are established in Section \ref{section_fixed_point_estimator}, with Section \ref{section_fixed_point_estimator_moderately_irregular_signal} devoted to the case $H \in (1/2,1)$, and Section \ref{section_fixed_point_estimator_rough_signal} addressing the case $H \in (1/3,1/2]$. Section \ref{section_implementation_simulations} presents some basic numerical experiments, while Section \ref{section_proofs} collects the proofs of our main results. Finally, to help the reader, all the key notations are summarized in Appendix \ref{key_notations_section}.
%


%
\section{Probabilistic preliminaries}\label{section_probabilistic_preliminaries}
As mentioned in the introduction, our estimation procedure for $\theta_0$ relies on tools from Malliavin calculus as well as Young/rough differential equations theory. For background on these topics, we refer the reader to \cite{NUALART06} and \cite{FH14} respectively.

This section builds on these classical results to develop several probabilistic statements required for our statistical analysis in Section \ref{section_fixed_point_estimator}. In particular, as already emphasized, a crucial step for our methodology is to reformulate the Malliavin derivative of $X_t$ ($t\in (0,T]$) in a form that does not explicitly depend on $B$, and to handle 2D Young integrals arising in the analysis.\\

In the sequel, we adopt the following standard non-degeneracy assumption on $\sigma$:
%


%
\begin{assumption}\label{non_degeneracy_condition}
The function $\sigma$ is bounded and $\inf_{\mathbb R}|\sigma| > 0$.
\end{assumption}
We are now ready to state the following proposition, which provides an expression for the Malliavin derivative of $X_t$ that is independent of $B$. This result will play a central role in establishing our theoretical guarantees on $\overline\theta_N$, and applies both for $H\in (1/3,1/2]$ and $H > 1/2$.
%


%
\begin{proposition}\label{Malliavin_derivative_X}
Under Assumption \ref{non_degeneracy_condition}, for every $s,t\in [0,T]$,
\begin{displaymath}
\mathbf D_sX_t =
\sigma(X_t)\exp\left(
\theta_0\int_{s}^{t}\left(b'(X_u) -
\frac{\sigma'(X_u)b(X_u)}{\sigma(X_u)}\right)du\right)\mathbf 1_{[0,t)}(s).
\end{displaymath}
\end{proposition}
This expression for the Malliavin derivative is crucial in transitioning from the fakestimator to the fixed-point estimator. 

We begin by considering the case $H \in (1/2,1)$, where the relationship between the Skorokhod and pathwise integrals is well established (see Proposition 5.2.3 in \cite{NUALART06}). Let us define
\begin{equation}\label{definition_phi_psi}
\pi = b\sigma,\quad
\varphi =\sigma\pi' =\sigma(\sigma b' +\sigma'b)\quad {\rm and}\quad
\psi =\frac{1}{\sigma}(\sigma b' -\sigma'b).   
\end{equation}
Then, combining this relationship with Proposition \ref{Malliavin_derivative_X} yields
\begin{eqnarray}
 \int_{0}^{T}b(X_s)\delta X_s
 & = &
 \theta_0\int_{0}^{T}b(X_s)^2ds +
 \int_{0}^{T}b(X_s)\sigma(X_s)\delta B_s
 \nonumber\\
 & = &
 \int_{0}^{T}b(X_s)dX_s -
 \alpha_H\int_{0}^{T}\int_{0}^{T}\mathbf D_s[\pi(X_t)]\cdot |t - s|^{2H - 2}dsdt
 \nonumber\\
 \label{relationship_Skorokhod_Young_X}
 & = &
 \int_{0}^{T}b(X_s)dX_s
 -\alpha_H\int_{0}^{T}\int_{0}^{t}\varphi(X_t)
 \exp\left(
 \theta_0\int_{s}^{t}\psi(X_u)du\right)|t - s|^{2H - 2}dsdt,
\end{eqnarray}
where $\alpha_H = H(2H - 1)$.

The situation is more challenging in the case $H\in (1/3,1/2]$. Here, it is again necessary to establish a link between the Skorokhod and pathwise integrals. In this context, Theorem 3.1 in \cite{ST22} is instrumental, as it provides an explicit relationship between the two, involving a correction term expressed through a two-dimensional integral.

In order to handle this correction term and to establish statistical guarantees for the fixed-point estimator, the following proposition will play a crucial role.
%


%
\begin{proposition}\label{2D_Young_integral_R}
Consider
\begin{displaymath}
\Delta_T =\{(s,t)\in [0,T]^2 : s < t\}
\quad\textrm{and}\quad
\widetilde\Delta_T =\{(s,t,u,v)\in [0,T]^4 : s < u < t < v\}.
\end{displaymath}
Assume that $H\in (1/3,1/2]$, and let $x : [0,T]^2\rightarrow\mathbb R$ be a function such that
\begin{equation}\label{2D_Young_integral_R_1}
|x(s,t)|\leqslant\mathfrak c_x|t - s|^{\alpha},
\quad\forall (s,t)\in\Delta_T,
\end{equation}
and
\begin{equation}\label{2D_Young_integral_R_2}
|x(s,t) - x(u,v)|\leqslant\mathfrak c_x(|s - u|^{\alpha} + |t - v|^{\alpha}),
\quad\forall (s,t,u,v)\in\widetilde\Delta_T,
\end{equation}
where $\alpha\in (1 - 2H,H)$ and $\mathfrak c_x$ is a positive constant. Then, the 2D Young integral of $x$ with respect to the covariance function $R$ of the fBm is well-defined, and
\begin{displaymath}
\int_{0 < s < t < T}x(s,t)dR(s,t) =
\alpha_H\int_{0}^{T}\int_{0}^{t}x(s,t)|t - s|^{2H - 2}dsdt.
\end{displaymath}
\end{proposition}
We have stated Proposition \ref{2D_Young_integral_R} here because we believe it deserves independent attention, as it may be applied in several contexts beyond our specific drift parameter estimation problem. In the following proposition, we instead detail its application within our framework. In particular, it provides an explicit expression for the Skorokhod integral of $b\circ X$ with respect to $X$ in the case $H\in (1/3,1/2]$. This result holds provided the drift function $b$ is bounded, as formalized below.
%


%
\begin{assumption}\label{bounded_drift_condition}
The function $b$ is bounded (i.e. $\mathfrak m_b :=\sup_{\mathbb R}|b| <\infty$).
\end{assumption}
%


%
\begin{proposition}\label{relationship_Skorokhod_rough}
Assume that $H\in (1/3,1/2]$. Under Assumptions \ref{non_degeneracy_condition} and \ref{bounded_drift_condition}, if $b$ is bounded, then
\begin{eqnarray}
 \label{relationship_Skorokhod_rough_1}
 \int_{0}^{T}\pi(X_s)\delta B_s
 & = &
 \int_{0}^{T}\pi(X_s)dB_s - H\int_{0}^{T}\varphi(X_s)s^{2H - 1}ds\\
 & &
 \hspace{2cm} -\alpha_H\int_{0}^{T}\int_{0}^{t}
 \varphi(X_t)\left(\exp\left(\theta_0\int_{s}^{t}\psi(X_u)du\right) - 1\right)|t - s|^{2H - 2}dsdt
 \nonumber
\end{eqnarray}
and
\begin{eqnarray}
 \label{relationship_Skorokhod_rough_2}
 \int_{0}^{T}b(X_s)\delta X_s
 & = &
 \int_{0}^{T}b(X_s)dX_s - H\int_{0}^{T}\varphi(X_s)s^{2H - 1}ds\\
 & &
 \hspace{2cm} -\alpha_H\int_{0}^{T}\int_{0}^{t}
 \varphi(X_t)\left(\exp\left(\theta_0\int_{s}^{t}\psi(X_u)du\right) - 1\right)|t - s|^{2H - 2}dsdt.
 \nonumber
\end{eqnarray}
\end{proposition}
%


%
\begin{remark}\label{remark_relationship_Skorokhod_rough}
Let us add a few remarks on Proposition \ref{relationship_Skorokhod_rough}.
\begin{enumerate}
 \item Since $\sigma$, $1/\sigma$, and the derivatives of $b$ and $\sigma$ are already assumed to be bounded, if $b$ is itself bounded, then so are $\varphi$ and $\psi$. This additional boundedness condition on $b$ (see Assumption \ref{bounded_drift_condition}), which is not needed when $H \in (1/2,1)$, is the price to pay in order to apply Theorem 3.1 of Song and Tindel as well as Proposition \ref{2D_Young_integral_R} in the proof of Proposition \ref{relationship_Skorokhod_rough} (see Section \ref{section_proofs}).
 \item It is worth noting that, after application of Proposition \ref{2D_Young_integral_R}, the relationships between Skorokhod and pathwise integrals in the two cases $H > 1/2$ and $H\in (1/3,1/2]$, summarized in Equations (\ref{relationship_Skorokhod_Young_X}) and (\ref{relationship_Skorokhod_rough_2}) respectively, appear quite similar. However, there is a crucial difference: in Equation (\ref{relationship_Skorokhod_rough_2}), the right-hand side cannot be simplified in such a way that the second integral cancels the $-1$ appearing in the third term, as happens in Equation (\ref{relationship_Skorokhod_Young_X}). This is because the map
 \begin{displaymath}
 (s,t)\in [0,T]^2\longmapsto
 \varphi(X_t)\exp\left(\theta_0\int_{s}^{t}\psi(X_u)du\right)\mathbf 1_{[0,t)}(s)
 \end{displaymath}
 does not satisfy the assumptions of Proposition \ref{2D_Young_integral_R}.
\end{enumerate}
\end{remark}
In the sequel, assume that $X^1,\dots,X^N$ have been observed on $[0,T_0]$ ($T_0 > 0$), and then on $[0,T]$ for any $T\in (0,T_0)$. Assume also that for every $t\in (0,T]$, the probability distribution of $X_t$ has a density $f_t$ with respect to the Lebesgue measure on $(\mathbb R,\mathcal B(\mathbb R))$ such that $s\mapsto f_s(x)$ ($x\in\mathbb R$) belongs to $\mathbb L^1([0,T],dt)$. Under Assumption \ref{non_degeneracy_condition}, this condition on the distribution of $X_t$ ($t\in (0,T]$) is satisfied when $\sigma$ is constant by \cite{LPS23}, Theorem 1.3, and when $\sigma$ is not constant but $b$ is bounded by \cite{BNOT16}, Theorem 1.5.

The primary purpose of this integrability condition on $s\mapsto f_s(x)$ ($x\in\mathbb R$) is to allow us to rigorously define the "time-averaged" density function $f$ given by
\begin{displaymath}
f(x) =\frac{1}{T}\int_{0}^{T}f_s(x)ds,
\quad\forall x\in\mathbb R.
\end{displaymath}
The usual norm on $\mathbb L^2(\mathbb R,f(x)dx)$ is denoted by $\|.\|_f$. Let us clarify that $f$ serves essentially as a convenient technical and notational device. It allows us to express asymptotic norms, such as $\|b\|_f$, and limiting variances in a compact and statistically interpretable manner.
%


%
\section{The fixed-point estimator of the drift parameter}\label{section_fixed_point_estimator}
By the (standard) law of large numbers, and observing that the Skorokhod integral is centered, we obtain
\begin{displaymath}
\widehat\theta_N =
\theta_0 +\left(\sum_{i = 1}^{N}\int_{0}^{T}b(X_{s}^{i})^2ds\right)^{-1}
\left(\sum_{i = 1}^{N}\int_{0}^{T}b(X_{s}^{i})\sigma(X_{s}^{i})\delta B_{s}^{i}\right)
\xrightarrow[N\rightarrow\infty]{\rm a.s.}\theta_0,
\end{displaymath}
so that $\widehat\theta_N$ is a consistent estimator of $\theta_0$. However, as already mentioned, $\widehat\theta_N$ cannot be computed directly from the data. Thanks to the relationships provided in Equations (\ref{relationship_Skorokhod_Young_X}) and (\ref{relationship_Skorokhod_rough_2}), we can replace the Skorokhod integral appearing in the fakestimator by the corrected pathwise integral, thereby moving to a computable approximation of $\widehat\theta_N$.

The goal of this section is to provide theoretical guarantees for such a computable approximation of $\widehat\theta_N$, defined as the fixed-point of a suitably chosen random functional. Due to Remark \ref{remark_relationship_Skorokhod_rough}.(2), the cases $H\in (1/2,1)$ and $H\in (1/3,1/2]$ need to be treated separately.
%


%
\subsection{The case $H\in (1/2,1)$}\label{section_fixed_point_estimator_moderately_irregular_signal}
First, let us define the following quantities:
\begin{displaymath}
D_N =\frac{1}{NT}\sum_{i = 1}^{N}\int_{0}^{T}b(X_{s}^{i})^2ds
\end{displaymath}
and
\begin{eqnarray}
 \label{pathwise_counterpart_I}
 I_N & = &
 \frac{1}{NTD_N}
 \sum_{i = 1}^{N}\int_{0}^{T}b(X_{s}^{i})dX_{s}^{i}\\
 & = &
 \frac{1}{NTD_N}
 \sum_{i = 1}^{N}(\texttt b(X_{T}^{i}) -\texttt b(x_0))
 \quad {\rm with}\quad\texttt b' = b,
 \nonumber
\end{eqnarray}
where the last equality follows from the change-of-variable formula for the Young integral. By Equality (\ref{relationship_Skorokhod_Young_X}), we can rewrite $\widehat\theta_N - I_N$ in the following way:
\begin{equation}\label{fakestimator_rewritten}
\widehat\theta_N - I_N =
-\frac{\alpha_H}{NTD_N}
\sum_{i = 1}^{N}\int_{0}^{T}\int_{0}^{t}
\varphi(X_{t}^{i})\exp\left(
\theta_0\int_{s}^{t}\psi(X_{u}^{i})du\right)
|t - s|^{2H - 2}dsdt.
\end{equation}
Consider the random functional
\begin{displaymath}
\Theta_N(\cdot) :=
-\frac{\alpha_H}{NTD_N}\sum_{i = 1}^{N}
\int_{0}^{T}\int_{0}^{t}\varphi(X_{t}^{i})
\exp\left((\cdot + I_N)\int_{s}^{t}\psi(X_{u}^{i})du\right)|t - s|^{2H - 2}dsdt.
\end{displaymath}
From Equality (\ref{fakestimator_rewritten}), and since $\widehat\theta_N$ is a consistent estimator of $\theta_0$,
\begin{displaymath}
\widehat\theta_N - I_N =
\Theta_N(\theta_0 - I_N)\approx
\Theta_N(\widehat\theta_N - I_N),
\end{displaymath}
which motivates introducing the estimator $\overline\theta_N = I_N + R_N$, where $R_N$ is (when it exists and is unique) a fixed-point of the map $\Theta_N$. In order to ensure the existence of a unique fixed-point, we impose the following "sign conditions".
%


%
\begin{assumption}\label{fixed_point_conditions}
The functions $b'$, $\varphi$ and $\psi$ are nonpositive, and $\theta_0 > 0$.
\end{assumption}
%


%
\begin{proposition}\label{existence_uniqueness_R_N_Young}
Under Assumptions \ref{non_degeneracy_condition} and \ref{fixed_point_conditions}, if
\begin{equation}\label{existence_uniqueness_R_N_Young_1}
T^{2H}\frac{M_N}{D_N}\leqslant
\frac{\mathfrak c}{\overline\alpha_H
\|\varphi\|_{\infty}\|\psi\|_{\infty}},
\end{equation}
where $\mathfrak c$ is a deterministic constant arbitrarily chosen in $(0,1)$,
\begin{displaymath}
M_N = e^{\|\psi\|_{\infty}|I_N|T}
\quad\textrm{and}\quad
\overline\alpha_H =\frac{|\alpha_H|}{2H(2H + 1)},
\end{displaymath}
then $\Theta_N$ is a contraction from $\mathbb R_+$ into itself. Therefore, $R_N$ exists and is unique.
\end{proposition}
%


%
\begin{remark}\label{remark_existence_uniqueness_R_N_Young}
Note that, in particular, the conditions of Proposition \ref{existence_uniqueness_R_N_Young} on $(\varphi,\psi)$ imply that it is bounded. Indeed,
\begin{displaymath}
\varphi =\sigma^2\left(b' +\frac{\sigma'}{\sigma}b\right)
\quad\textrm{and}\quad
\psi = b' -\frac{\sigma'}{\sigma}b,
\end{displaymath}
leading to
\begin{displaymath}
\left\{
\begin{array}{rcl}
 b' & \leqslant & 0\\
 \varphi & \leqslant & 0\\
 \psi & \leqslant & 0
\end{array}\right.
\Longleftrightarrow
\left\{
\begin{array}{rcl}
 b' & \leqslant & 0\\
 -b' -\sigma'b/\sigma & \geqslant & 0\\
 b' -\sigma'b/\sigma & \leqslant & 0
\end{array}\right.
\Longleftrightarrow
\left\{
\begin{array}{rcl}
 b' & \leqslant & 0\\
 b' & \leqslant & \sigma'b/\sigma\leqslant -b'
\end{array}\right.
\quad ({\rm A})
\end{displaymath}
Since $b'$ and $\sigma$ are bounded, the ratio $\sigma'b/\sigma$ is also bounded, and therefore so are $(\varphi,\psi)$.
\end{remark}
%


%
\begin{example}\label{example_existence_uniqueness_R_N_Young}
Let us provide examples of drift and volatility functions satisfying the condition ${\rm (A)}$:
\begin{enumerate}
 \item If $\sigma$ is constant, then $(\varphi,\psi) = (\sigma^2b',b')$ satisfies ${\rm (A)}$ if and only if $b'\leqslant 0$ (as in \cite{MARIE25a}, Proposition 7).
 \item Assume that $b =\sigma = {\tt F}$, where ${\tt F}\in C^1(\mathbb R)$ satisfies
 \begin{displaymath}
 \inf_{x\in\mathbb R}|{\tt F}(x)| > 0,\quad
 \sup_{x\in\mathbb R}|{\tt F}(x)| <\infty,\quad
 {\tt F}'\leqslant 0
 \quad {\rm and}\quad
 \sup_{x\in\mathbb R}|{\tt F}'(x)| <\infty.
 \end{displaymath}
 Then, $b'\leqslant 0$ and
 \begin{displaymath}
 \frac{\sigma'b}{\sigma} = {\tt F}' = b'.
 \end{displaymath}
 So, $(\varphi,\psi)$ satisfies the condition (A), and Equation (\ref{main_equation}) is equivalent to
 \begin{displaymath}
 X_t = x_0 +\int_{0}^{t}{\tt F}(X_s)(\theta_0ds + dB_s),
 \quad t\in [0,T].
 \end{displaymath}
 \item Assume that $b(x) = -x$, which is quite common. Then,
 \begin{displaymath}
 ({\rm A})\Longleftrightarrow
 -1\leqslant\frac{\sigma'(x)}{\sigma(x)}x\leqslant 1.
 \end{displaymath}
 For instance,
 \begin{itemize}
  \item If $\sigma(x) =\pi +\arctan(x)$, then $\sigma'(x) = (1 + x^2)^{-1}$, leading to
  \begin{displaymath}
  \frac{\sigma'(x)}{\sigma(x)}x =
  \frac{x}{(1 + x^2)(\pi +\arctan(x))}\in [-1,1].
  \end{displaymath}
  \item If $\sigma(x) = 1 + e^{-x^2}$, then $\sigma'(x) = -2xe^{-x^2}$, leading to
  \begin{displaymath}
  \frac{\sigma'(x)}{\sigma(x)}x =\frac{-2x^2e^{-x^2}}{1 + e^{-x^2}}\in [-1,1].
  \end{displaymath}
 \end{itemize}
\end{enumerate}
\end{example}
The previous examples illustrate that, although the conditions imposed on the coefficients for our analysis may appear rather restrictive, there is nonetheless a broad class of drift and volatility functions satisfying them.

As established in Proposition \ref{existence_uniqueness_R_N_Young}, Assumption \ref{fixed_point_conditions} and the condition (\ref{existence_uniqueness_R_N_Young_1}) are sufficient to ensure the existence and uniqueness of the fixed-point $R_N$. Consequently, we maintain Assumption \ref{fixed_point_conditions} in the sequel, and restrict our study to the event where (\ref{existence_uniqueness_R_N_Young_1}) is satisfied. To this end, let us define
\begin{equation}\label{Delta_event}
\Delta_N :=\left\{T^{2H}\frac{M_N}{D_N}\leqslant
\frac{\mathfrak c}{\overline\alpha_H
\|\varphi\|_{\infty}\|\psi\|_{\infty}}\right\}.
\end{equation}
We thus introduce a modified estimator, $\overline\theta_{N}^{\mathfrak c} :=\overline\theta_N\mathbf 1_{\Delta_N}$, which coincides with the original estimator on the set where the fixed-point of $\Theta_N$ is well-defined. To ensure that $\overline\theta_{N}^{\mathfrak c}$ inherits the asymptotic properties of $\widehat\theta_N$, it remains to show that the probability of the complement event $\mathbb P(\Delta_{N}^{c})$ vanishes at an appropriate rate.
%


%
\begin{proposition}\label{deviation_bound}
Assume that
\begin{equation}\label{deviation_bound_1}
\frac{2T^{2H}}{\|b\|_{f}^{2}}\exp\left(\frac{2\|\psi\|_{\infty}}{\|b\|_{f}^{2}}
|\mathbb E({\tt b}(X_T)) - {\tt b}(x_0)|\right)
<\frac{\mathfrak c}{\overline\alpha_H
\|\varphi\|_{\infty}\|\psi\|_{\infty}}.
\end{equation}
Under Assumptions \ref{non_degeneracy_condition} and \ref{fixed_point_conditions}, there exists a constant $\mathfrak c_{\ref{deviation_bound}} > 0$, not depending on $N$, such that
\begin{displaymath}
\mathbb P(\Delta_{N}^{c})
\leqslant
\frac{\mathfrak c_{\ref{deviation_bound}}}{N}.
\end{displaymath}
\end{proposition}
Now, the following proposition provides an asymptotic confidence interval for the computable fixed-point estimator $\overline\theta_{N}^{\mathfrak c}$.
\\
\\
{\bf Notation.} $u_{\cdot} =\phi^{-1}(\cdot)$, where $\phi$ is the standard normal distribution function.
%


%
\begin{proposition}\label{ACI_approximation_LS_Young}
Under Assumptions \ref{non_degeneracy_condition} and \ref{fixed_point_conditions}, for any $\alpha\in (0,1/2)$, if $T$ satisfies the condition (\ref{deviation_bound_1} with
\begin{displaymath}
\mathfrak c\in\left(0,
1 -\frac{u_{1 -\frac{\lambda\alpha}{2}}}{u_{1 -\frac{\alpha}{2}}}\right)
\quad\textrm{and}\quad
\lambda\in\left(1,\frac{1}{2\alpha}\right),
\end{displaymath}
then
\begin{displaymath}
\lim_{N\rightarrow\infty}
\mathbb P\left(\theta_0\in\left[
\overline\theta_{N}^{\mathfrak c} -
\sqrt{\frac{Y_N}{ND_{N}^{2}}}u_{1 -\frac{\alpha}{2}}\textrm{ $;$ }
\overline\theta_{N}^{\mathfrak c} +
\sqrt{\frac{Y_N}{ND_{N}^{2}}}u_{1 -\frac{\alpha}{2}}\right]\right)
\geqslant 1 - 2\lambda\alpha,
\end{displaymath}
where
\begin{eqnarray*}
 Y_N & := &
 \frac{1}{NT^2}
 \sum_{i = 1}^{N}\left(
 \alpha_H\int_{0}^{T}\int_{0}^{T}
 \pi(X_{s}^{i})\pi(X_{t}^{i})|t - s|^{2H - 2}dsdt\right.\\
 & &
 \hspace{3.5cm}\left.
 +\alpha_{H}^{2}\int_{[0,T]^2}\int_{0}^{v}\int_{0}^{u}
 |u -\overline u|^{2H - 2}|v -\overline v|^{2H - 2}
 \varphi(X_{v}^{i})\varphi(X_{u}^{i})d\overline ud\overline vdudv\right).
\end{eqnarray*}
\end{proposition}
%


%
\begin{remark}\label{remark_ACI_approximation_LS_Young}
Consider $\alpha\in (0,1/2)$. As established in the proof of Proposition \ref{ACI_approximation_LS_Young}, a non-computable sharp asymptotic $(1 -\alpha)$-confidence interval for $\theta_0$ is given - for $N$ large enough - by
\begin{displaymath}
\mathbb I_{N}^{*} =\left[
\widehat\theta_N -
\sqrt{\frac{Y_{N}^{*}}{ND_{N}^{2}}}u_{1 -\frac{\alpha}{2}}\textrm{ $;$ }
\widehat\theta_N +
\sqrt{\frac{Y_{N}^{*}}{ND_{N}^{2}}}u_{1 -\frac{\alpha}{2}}\right],
\end{displaymath}
where
\begin{displaymath}
Y_{N}^{*} =
\frac{1}{NT^2}
\sum_{i = 1}^{N}\left(
\alpha_H\int_{0}^{T}\int_{0}^{T}
\pi(X_{s}^{i})\pi(X_{t}^{i})|t - s|^{2H - 2}dsdt + R_i\right)
\end{displaymath}
and, for every $i\in\{1,\dots,N\}$,
\begin{eqnarray*}
 R_i & = &
 \alpha_{H}^{2}\int_{[0,T]^2}\int_{0}^{v}\int_{0}^{u}
 |u -\overline u|^{2H - 2}|v -\overline v|^{2H - 2}\\
 & &
 \hspace{3cm}\times
 \varphi(X_{v}^{i})\varphi(X_{u}^{i})\exp\left(\theta_0\left(\int_{\overline u}^{v}\psi(X_{s}^{i})ds +
 \int_{\overline v}^{u}\psi(X_{s}^{i})ds\right)\right)
 d\overline ud\overline vdudv.
\end{eqnarray*}
The computable asymptotic $(1 - 2\lambda\alpha)$-confidence interval
\begin{displaymath}
\mathbb I_N =\left[
\overline\theta_{N}^{\mathfrak c} -
\sqrt{\frac{Y_N}{ND_{N}^{2}}}u_{1 -\frac{\alpha}{2}}\textrm{ $;$ }
\overline\theta_{N}^{\mathfrak c} +
\sqrt{\frac{Y_N}{ND_{N}^{2}}}u_{1 -\frac{\alpha}{2}}\right],
\end{displaymath}
provided in Proposition \ref{ACI_approximation_LS_Young}, is not as sharp as $\mathbb I_{N}^{*}$ but remains satisfactory:
\begin{itemize}
 \item Since $\lambda$ can be chosen as close to $1$ as possible, and since $\alpha$ is near of $0$ in practice, the confidence level of $\mathbb I_N$ is not that degraded with respect to that of $\mathbb I_{N}^{*}$. Indeed, for $\alpha = 1\%$ (resp. $\alpha = 5\%$), $\theta_0$ belongs to $\mathbb I_{N}^{*}$ with probability $0.99$ (resp. $0.95$), while $\theta_0$ belongs to $\mathbb I_N$ with a probability of around $0.98$ (resp. $0.90$). 
 \item Let us compare the radii
 \begin{displaymath}
 {\tt R}(\mathbb I_N) =
 \sqrt{\frac{Y_N}{ND_{N}^{2}}}u_{1 -\frac{\alpha}{2}}
 \quad\textrm{and}\quad
 {\tt R}(\mathbb I_{N}^{*}) =
 \sqrt{\frac{Y_{N}^{*}}{ND_{N}^{2}}}u_{1 -\frac{\alpha}{2}}.
 \end{displaymath}
 By the definition of $Y_N$ and $Y_{N}^{*}$,
 \begin{displaymath}
 |{\tt R}(\mathbb I_N) - {\tt R}(\mathbb I_{N}^{*})|\leqslant
 \frac{u_{1 -\frac{\alpha}{2}}}{\sqrt ND_N}\Delta_N(\theta_0)^{\frac{1}{2}}
 \end{displaymath}
 with
 \begin{eqnarray*}
  \Delta_N(\theta_0) & := &
  \frac{\alpha_{H}^{2}}{NT^2}
  \sum_{i = 1}^{N}
  \int_{[0,T]^2}\int_{0}^{v}\int_{0}^{u}
  |u -\overline u|^{2H - 2}|v -\overline v|^{2H - 2}\\
  & &
  \hspace{0.75cm}\times
  \varphi(X_{v}^{i})\varphi(X_{u}^{i})
  \left|1 -\exp\left(\theta_0\left(\int_{\overline u}^{v}\psi(X_{s}^{i})ds +
  \int_{\overline v}^{u}\psi(X_{s}^{i})ds\right)\right)\right|
  d\overline ud\overline vdudv\\
  & \leqslant &
  \frac{\alpha_{H}^{2}}{NT^2}
  \sum_{i = 1}^{N}
  \int_{[0,T]^2}\int_{0}^{v}\int_{0}^{u}
  |u -\overline u|^{2H - 2}|v -\overline v|^{2H - 2}
  \varphi(X_{v}^{i})\varphi(X_{u}^{i})
  d\overline ud\overline vdudv.
 \end{eqnarray*}
 Moreover, since $\varphi$ is bounded under the conditions of Proposition \ref{ACI_approximation_LS_Young} (see Remark \ref{remark_existence_uniqueness_R_N_Young}),
 \begin{displaymath}
 \Delta_N(\theta_0)\leqslant
 \frac{\|\varphi\|_{\infty}^{2}}{T^2}
 \left(\alpha_H\int_{0}^{T}\int_{0}^{u}
 |u -\overline u|^{2H - 2}d\overline udu\right)^2\leqslant
 \|\varphi\|_{\infty}^{2}T^{2(2H - 1)},
 \end{displaymath}
 leading to
 \begin{displaymath}
 |{\tt R}(\mathbb I_N) - {\tt R}(\mathbb I_{N}^{*})|\leqslant
 \frac{u_{1 -\frac{\alpha}{2}}\|\varphi\|_{\infty}^{2}}{\sqrt ND_N}T^{2(2H - 1)}.
 \end{displaymath}
 Regardless of the values of $N$ and $\alpha$, since $T^{2(2H - 1)}\downarrow 0$ when $T$ does, the radius of $\mathbb I_N$ is close to that of the sharp $(1 -\alpha)$-confidence interval $\mathbb I_{N}^{*}$ by taking a small time horizon $T$, which is already required due to the condition (\ref{deviation_bound_1}).
\end{itemize}
\end{remark}
Finally, in order to establish a non-asymptotic risk bound on our estimator, we must bound its denominator away from zero. To this end, we introduce the following truncated estimator:
\begin{equation}\label{fixed_point_truncated_estimator}
\overline\theta_{N}^{\mathfrak c,\mathfrak d} =
\overline\theta_{N}^{\mathfrak c}\mathbf 1_{D_N\geqslant\mathfrak d}
\quad {\rm with}\quad
\mathfrak d\in\Delta_f :=\left(0,\frac{\|b\|_{f}^{2}}{2}\right].
\end{equation}
%


%
\begin{proposition}\label{risk_bound_approximation_LS_Young}
Under Assumptions \ref{non_degeneracy_condition} and \ref{fixed_point_conditions}, if $T$ satisfies the condition (\ref{deviation_bound_1}), then there exists a constant $\mathfrak c_{\ref{risk_bound_approximation_LS_Young}} > 0$, not depending on $N$, such that
\begin{displaymath}
\mathbb E(|\overline\theta_{N}^{\mathfrak c,\mathfrak d} -\theta_0|^2)\leqslant
\frac{\mathfrak c_{\ref{risk_bound_approximation_LS_Young}}}{N}.
\end{displaymath}
\end{proposition}
In this way, in the case $H > 1/2$, we have established both an asymptotic confidence interval and a non-asymptotic risk bound for the fixed-point estimator. Let us now turn to the case $H\in (1/3,1/2]$.
%


%
\subsection{The case $H\in (1/3,1/2]$}\label{section_fixed_point_estimator_rough_signal}
First, as in the case $H > 1/2$, we start by defining the functional underlying our fixed-point estimator. Note that, by applying the change-of-variable formula for the rough integral (instead of the Young integral), the quantity $I_N$ retains the same form (\ref{pathwise_counterpart_I}) as in the case $H > 1/2$. However, the connection between the Skorokhod and pathwise integrals, previously expressed in Equality (\ref{relationship_Skorokhod_Young_X}), is now replaced by the relation in Equality (\ref{relationship_Skorokhod_rough_2}). This leads us to introduce a modified version of the random functional $\Theta_N$ from the previous section:
\begin{displaymath}
\widetilde\Theta_N(\cdot) :=
\frac{1}{NTD_N}\sum_{i = 1}^{N}
\int_{0}^{T}\varphi(X_{t}^{i})\underbrace{\left[-Ht^{2H - 1} +
\alpha_H\int_{0}^{t}\left(1 -\exp\left((\cdot + I_N)
\int_{s}^{t}\psi(X_{u}^{i})du\right)\right)|t - s|^{2H - 2}ds\right]}_{=:\Lambda_{t}^{i}(\cdot + I_N)}dt.
\end{displaymath}
Again, since $\widehat\theta_N$ is a consistent estimator of $\theta_0$, we have
\begin{displaymath}
\widehat\theta_N - I_N =
\widetilde\Theta_N(\theta_0 - I_N)\approx
\widetilde\Theta_N(\widehat\theta_N - I_N),
\end{displaymath}
which justifies considering the estimator $\overline\theta_N = I_N + R_N$ of $\theta_0$, where $R_N$ is the fixed-point (when it exists and is unique) of the map $\widetilde\Theta_N$.
%


%
\begin{proposition}\label{existence_uniqueness_R_N_rough}
Under Assumptions \ref{non_degeneracy_condition}, \ref{bounded_drift_condition} and \ref{fixed_point_conditions}, if $M_N/D_N$ satisfies the condition (\ref{existence_uniqueness_R_N_Young_1}), then $\widetilde\Theta_N$ is a contraction from $\mathbb R_+$ into itself. Therefore, $R_N$ exists and is unique.
\end{proposition}
As in the case $H > 1/2$, the existence and uniqueness of the fixed-point of the map $\widetilde\Theta_N$ are crucial to establish theoretical guarantees on our estimator of $\theta_0$. Consequently, we maintain Assumptions \ref{non_degeneracy_condition}, \ref{bounded_drift_condition} and \ref{fixed_point_conditions} in the sequel. We now shift our focus to the estimator $\overline\theta_{N}^{\mathfrak c} =\overline\theta_N\mathbf 1_{\Delta_N}$, where the event $\Delta_N$ remains defined as in Equality (\ref{Delta_event}). On this event, the condition (\ref{existence_uniqueness_R_N_Young_1}) is automatically satisfied. Again, to transfer suitable properties from $\widehat\theta_N$ to $\overline\theta_{N}^{\mathfrak c}$, it is necessary to control $\mathbb P(\Delta_{N}^{c})$. This control is already established in Proposition \ref{deviation_bound}, and a brief review of its proof confirms that the result does not depend on $H$ nor on the considered functional ($\Theta_N$ or $\widetilde\Theta_N$).

Now, we proceed to the following proposition providing an asymptotic confidence interval for the computable fixed-point estimator $\overline\theta_{N}^{\mathfrak c}$.
%


%
\begin{proposition}\label{ACI_approximation_LS_rough}
Assume that $\theta_0\in (0,\theta_{\max}]$ with a known $\theta_{\max} > 0$. Under Assumptions \ref{non_degeneracy_condition}, \ref{bounded_drift_condition} and \ref{fixed_point_conditions}, for any $\alpha\in (0,1/2)$, if $T$ satisfies the condition (\ref{deviation_bound_1}) with
\begin{displaymath}
\mathfrak c\in\left(0,
1 -\frac{u_{1 -\frac{\lambda\alpha}{2}}}{u_{1 -\frac{\alpha}{2}}}\right)
\quad\textrm{and}\quad
\lambda\in\left(1,\frac{1}{2\alpha}\right),
\end{displaymath}
then
\begin{displaymath}
\lim_{N\rightarrow\infty}
\mathbb P\left(\theta_0\in\left[
\overline\theta_{N}^{\mathfrak c} -
\sqrt{\frac{\mathfrak Y_N}{ND_{N}^{2}}}u_{1 -\frac{\alpha}{2}}\textrm{ $;$ }
\overline\theta_{N}^{\mathfrak c} +
\sqrt{\frac{\mathfrak Y_N}{ND_{N}^{2}}}u_{1 -\frac{\alpha}{2}}\right]\right)
\geqslant 1 - 2\lambda\alpha,
\end{displaymath}
where
\begin{displaymath}
\mathfrak Y_N =
\frac{1}{NT^2}\sum_{i = 1}^{N}\left(
\left|\int_{0}^{T}b(X_{s}^{i})dX_{s}^{i}\right| +
\theta_{\max}\int_{0}^{T}b(X_{s}^{i})^2ds +
\int_{0}^{T}\varphi(X_{t}^{i})\Lambda_{t}^{i}(\theta_{\max})dt
\right)^2.
\end{displaymath}
\end{proposition}
%


%
\begin{remark}\label{remark_ACI_approximation_LS_rough}
Let us make a few comments about Proposition \ref{ACI_approximation_LS_rough}.
\begin{enumerate}
 \item Comparing the (asymptotic) confidence interval in Proposition \ref{ACI_approximation_LS_Young} with that of Proposition \ref{ACI_approximation_LS_rough}, one can notice that $Y_N$ is now replaced by $\mathfrak Y_N$. This arises because, when $H > 1/2$, the variance of the Skorokhod integral with respect to $B$ can be directly calculated thanks to \cite{BHOZ08}, Theorem 3.11.1, which is no longer applicable when $H\in [1/3,1/2)$. So, $\mathfrak Y_N$ is obtained thanks to the relationship between the Skorokhod and the pathwise integrals provided in Equality (\ref{relationship_Skorokhod_rough_2}), leading to a more conservative confidence interval than that provided in Proposition \ref{ACI_approximation_LS_Young}.
 \item Note that $\theta_{\max}$ may be known even though $\theta_0$ is unknown; typically when $\theta_0$ is assumed to be a proportion (i.e. $\theta_0\in (0,\theta_{\max}]$ with $\theta_{\max} = 1$). For instance, if $b(x) = -x$, and if Equation (\ref{main_equation}) models the concentration of a drug administered intravenously to a patient involved in a clinical trial, then $\theta_0$ - the elimination rate constant - belongs to $(0,1]$. Otherwise, since $\widetilde\Theta_N$ doesn't depend on $\theta_{\max}$, then one can replace $\theta_{\max}$ by $\overline\theta_{N}^{\mathfrak c}$ - which is a converging estimator of $\theta_0$ - in the definition of $\mathfrak Y_N$ in practice.
\end{enumerate}
\end{remark}
Finally, the following proposition provides a non-asymptotic risk bound on the truncated estimator $\overline\theta_{N}^{\mathfrak c,\mathfrak d}$ defined as in (\ref{fixed_point_truncated_estimator}).
%


%
\begin{proposition}\label{risk_bound_approximation_LS_rough}
Under Assumptions \ref{non_degeneracy_condition}, \ref{bounded_drift_condition} and \ref{fixed_point_conditions}, if $T$ satisfies the condition (\ref{deviation_bound_1}), then there exists a constant $\mathfrak c_{\ref{risk_bound_approximation_LS_rough}} > 0$, not depending on $N$, such that
\begin{displaymath}
\mathbb E(|\overline\theta_{N}^{\mathfrak c,\mathfrak d} -\theta_0|^2)\leqslant
\frac{\mathfrak c_{\ref{risk_bound_approximation_LS_rough}}}{N}.
\end{displaymath}
\end{proposition}
An experienced reader might have anticipated that the proofs of Propositions \ref{risk_bound_approximation_LS_Young} and \ref{risk_bound_approximation_LS_rough} rely on the risk bounds for the fakestimator, together with a control on the approximation error incurred when moving from the fakestimator to the fixed-point estimator. The first part does not depend on $H$, while the second uses the bound on $\mathbb P(\Delta_{N}^{c})$ which, as discussed above, is also independent of $H$. It then follows directly that the proof of Proposition \ref{risk_bound_approximation_LS_rough} proceeds along the same lines as that of Proposition \ref{risk_bound_approximation_LS_Young}, and is therefore omitted.
%


%
\section{Implementation and simulations}\label{section_implementation_simulations}
%


%
\subsection{Implementation}\label{section_implementation}
Consider $n\in\mathbb N\backslash\{0,1\}$,
\begin{displaymath}
t_k =\frac{kT}{n},
\quad\forall k\in\{0,\dots,n\},
\end{displaymath}
and
\begin{displaymath}
t_{k,k + 1} = t_{k + 1} - t_k,
\quad\forall k\in\{0,\dots,n - 1\}.
\end{displaymath}
In practice, $\overline\theta_N$ may be approximated by
\begin{displaymath}
\overline\theta_{N,n,\overline n} = I_{N,n} + R_{N,n,\overline n},
\end{displaymath}
where
\begin{displaymath}
I_{N,n} =
\frac{1}{NTD_{N,n}}\sum_{i = 1}^{N}(\texttt b(X_{T}^{i}) -\texttt b(x_0))
\quad {\rm with}\quad
D_{N,n} =
\frac{1}{NT}
\sum_{i = 1}^{N}\sum_{k = 0}^{n - 1}b(X_{t_k}^{i})^2(t_{k + 1} - t_k),
\end{displaymath}
and
\begin{displaymath}
R_{N,n,\overline n} =
\left\{
\begin{array}{rcl}
 \underbrace{(\Theta_{N,n}\circ\dots\circ\Theta_{N,n})}_{\overline n\textrm{ times}}(0)
 & {\rm if} & H\in (1/2,1)\\
 \underbrace{
 (\widetilde\Theta_{N,n}\circ\dots\circ
 \widetilde\Theta_{N,n})}_{\overline n\textrm{ times}}(0)
 & {\rm if} & H\in (1/3,1/2]
\end{array}
\right.
\end{displaymath}
with $\overline n\in\mathbb N^*$,
\begin{displaymath}
\Theta_{N,n}(\cdot) :=
-\frac{\alpha_H}{NTD_{N,n}}
\sum_{i = 1}^{N}\sum_{k = 1}^{n - 1}\sum_{\ell = 0}^{k - 1}
\varphi(X_{t_k}^{i})
\exp\left((\cdot + I_{N,n})\sum_{r =\ell}^{k}\psi(X_{t_r}^{i})t_{r,r + 1}\right)
|t_k - t_{\ell}|^{2H - 2}t_{k,k + 1}t_{\ell,\ell + 1}
\end{displaymath}
and
\begin{eqnarray*}
 \widetilde\Theta_{N,n}(\cdot)
 & := &
 -\frac{1}{NTD_{N,n}}
 \sum_{i = 1}^{N}\sum_{k = 1}^{n - 1}
 \varphi(X_{t_k}^{i})\\
 & &
 \hspace{0.25cm}\times\left[Ht_{k}^{2H - 1} -
 \alpha_H\sum_{\ell = 0}^{k - 1}\left(1 -\exp\left((\cdot + I_{N,n})
 \sum_{r =\ell}^{k}\psi(X_{t_r}^{i})t_{r,r + 1}\right)\right)
 |t_k - t_{\ell}|^{2H - 2}t_{\ell,\ell + 1}\right]t_{k,k + 1}.
\end{eqnarray*}
%


%
\begin{remark}\label{remark_approximations}
All the integrals involved in the d\'efinitions of $D_N$, $I_N$, $\Theta_N$ and $\widetilde\Theta_N$ are replaced by the corresponding Riemann sums in those of $D_{N,n}$, $I_{N,n}$, $\Theta_{N,n}$ and $\widetilde\Theta_{N,n}$. Moreover, to approximate $R_N$ by $R_{N,n,\overline n}$ makes sense because when $H\in (1/2,1)$ (resp. $H\in (1/3,1/2]$), under the conditions of Proposition \ref{existence_uniqueness_R_N_Young} (resp. Proposition \ref{existence_uniqueness_R_N_rough}),
\begin{displaymath}
\underbrace{(\Theta_N\circ\dots\circ\Theta_N)}_{\overline n\textrm{ times}}(0)
\xrightarrow[\overline n\rightarrow\infty]{}R_N
\quad\textrm{(resp.}\quad
\underbrace{(\widetilde\Theta_N\circ\dots\circ\widetilde\Theta_N)}_{\overline n\textrm{ times}}(0)
\xrightarrow[\overline n\rightarrow\infty]{}R_N\textrm{)}
\end{displaymath}
thanks to Picard's fixed-point theorem.
\end{remark}
%


%
\subsection{Numerical experiments}\label{section_numerical_experiments}
In this section, some numerical experiments on the computable approximation of the least squares estimator of $\theta_0$ are presented for the three following models:
\begin{itemize}
 \item[(A)]\quad $dX_t = -\theta_0X_tdt + dB_t$ with $\theta_0 = 1$,
 \item[(B)]\quad $dX_t = -\theta_0X_tdt +\sigma_1(X_t)dB_t$ with $\theta_0 = 1$ and $\sigma_1(x) = 1 + e^{-x^2}$,
 \item[(C)]\quad $dX_t = -\theta_0X_tdt +\sigma_2(X_t)dB_t$ with $\theta_0 = 1$ and $\sigma_2(x) =\pi +\arctan(x)$.
\end{itemize}
For each model, with $H = 0.7$ and $H = 0.9$, $\overline\theta_N$ is computed from $N = 1,\dots,50$ paths of the process $X$. This experiment is repeated $100$ times. The means and the standard deviations of the error $|\overline\theta_{50} -\theta_0|$ are stored in Table \ref{table_errors}.
\begin{table}[h!]
\begin{center}
\begin{tabular}{lccc}
 \hline
  & $H$ & Mean error & Error StD.\\
 \hline
 Model (A) & 0.7 & 0.0489 & 0.0336\\
  & 0.9 & 0.0186 & 0.0139\\
 \hline
 Modal (B) & 0.7 & 0.0510 & 0.0500\\
  & 0.9 & 0.0294 & 0.0214\\
 \hline
 Model (C) & 0.7 & 0.1597 & 0.1443\\
  & 0.9 & 0.0823 & 0.0922\\
 \hline
\end{tabular}
\medskip
\caption{Means and standard deviations of the error of $\overline\theta_{50}$ (100 repetitions).}
\label{table_errors}
\end{center}
\end{table}
As expected, both the mean error and the standard deviation of the estimator $\overline{\theta}_{50}$ are larger for the models driven by multiplicative noise (Models (B) and (C)) compared to the fractional Ornstein-Uhlenbeck process (Model (A)). However, the mean error of $\overline\theta_{50}$ remains small for both Model (A) and Model (B): lower than $5.1\cdot 10^{-2}$. Since - for these two models - the error standard deviation of $\overline\theta_{50}$ is also small ($< 5\cdot 10^{-2}$), on average, the error of $\overline\theta_{50}$ for one repetition of the experiment should be near of its mean error. Finally, note also that for the three models ((A), (B) and (C)), the mean error of $\overline\theta_{50}$ is higher when $H = 0.7$ than when $H = 0.9$; probably because $H$ controls the H\"older exponent of the paths of the fractional Brownian motion.
\begin{figure}[h!]
\begin{center}
\begin{tabular}{cc}
\includegraphics[scale=0.25]{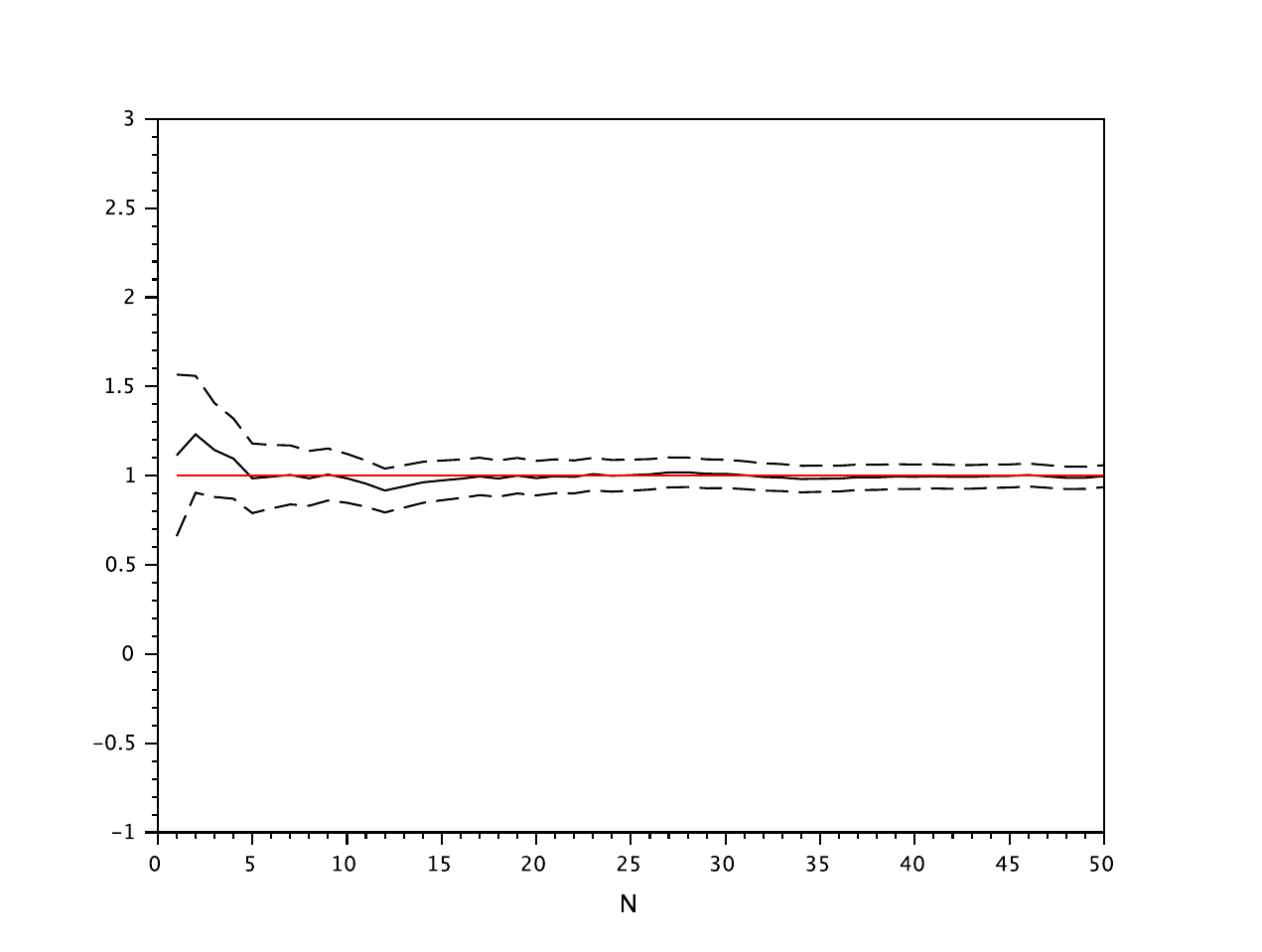} &
\includegraphics[scale=0.25]{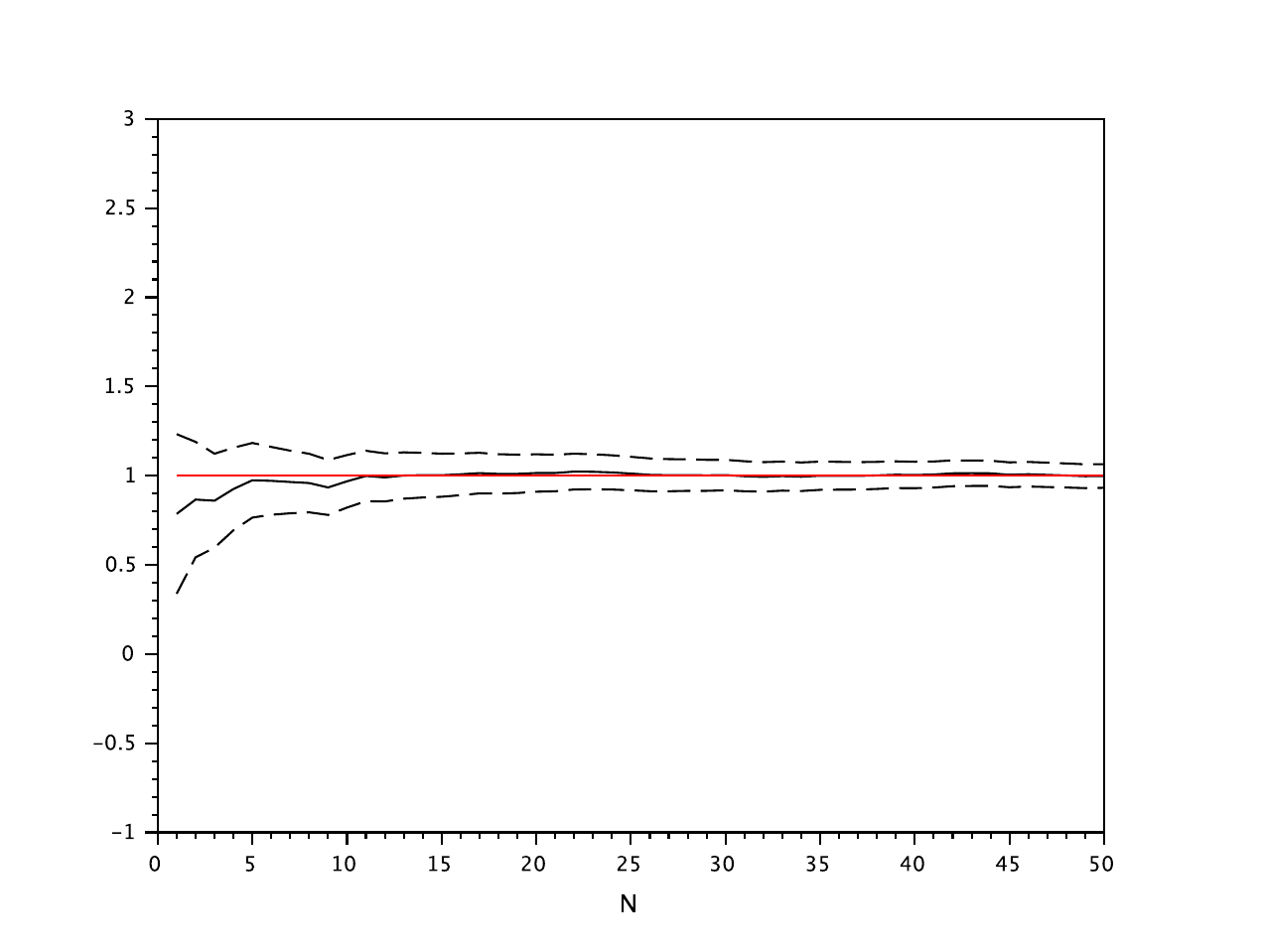} 
\end{tabular}
\end{center}
\caption{Plots of $N\mapsto\overline\theta_N$ (black line) and of the bounds of the $89\%$-ACIs (dashed black lines) for Model (A) with $H = 0.7$ (left) and $H = 0.9$ (right).}
\label{plots_Model_A}
\end{figure}
\begin{figure}[h!]
\begin{center}
\begin{tabular}{cc}
\includegraphics[scale=0.25]{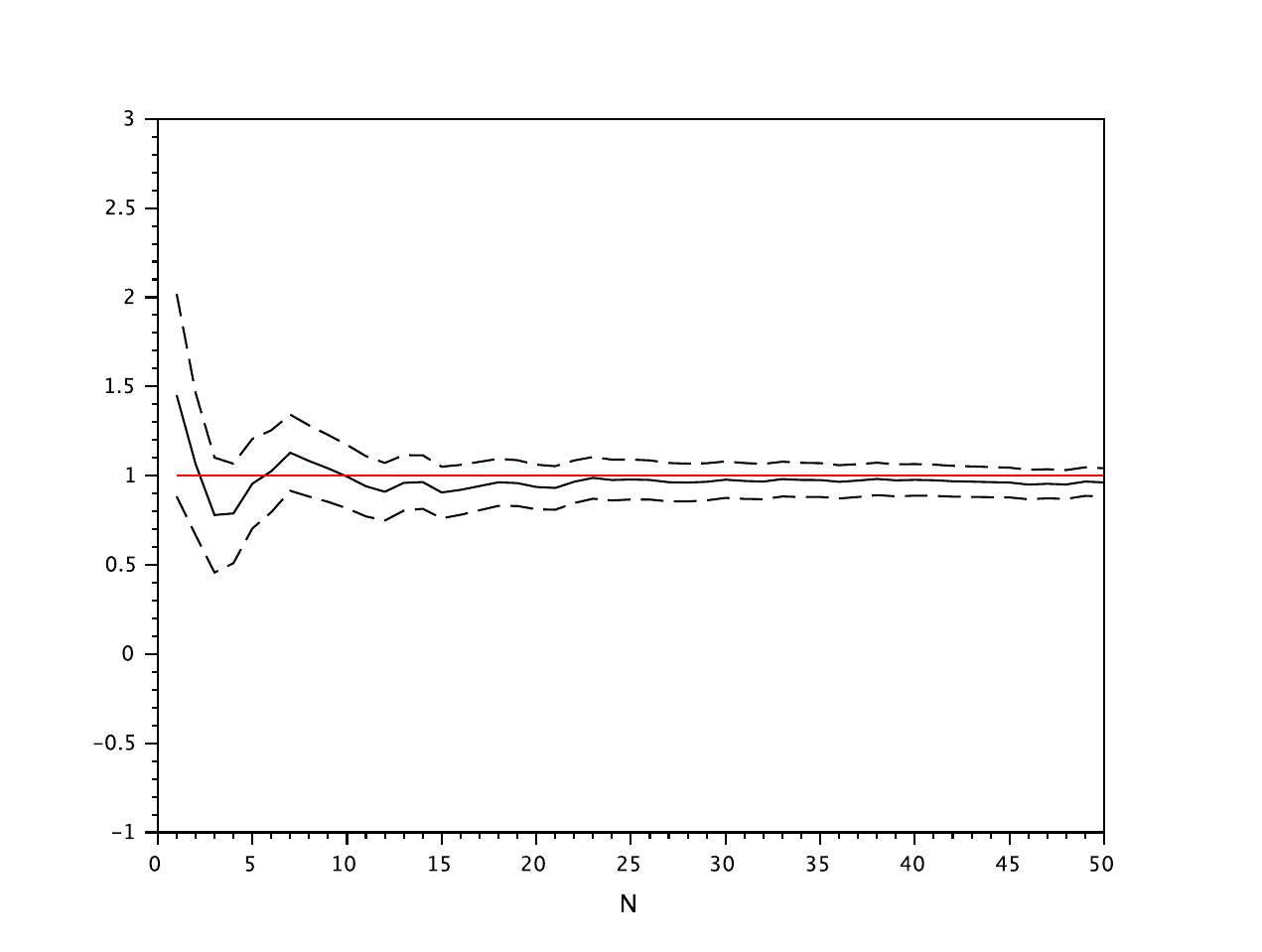} &
\includegraphics[scale=0.25]{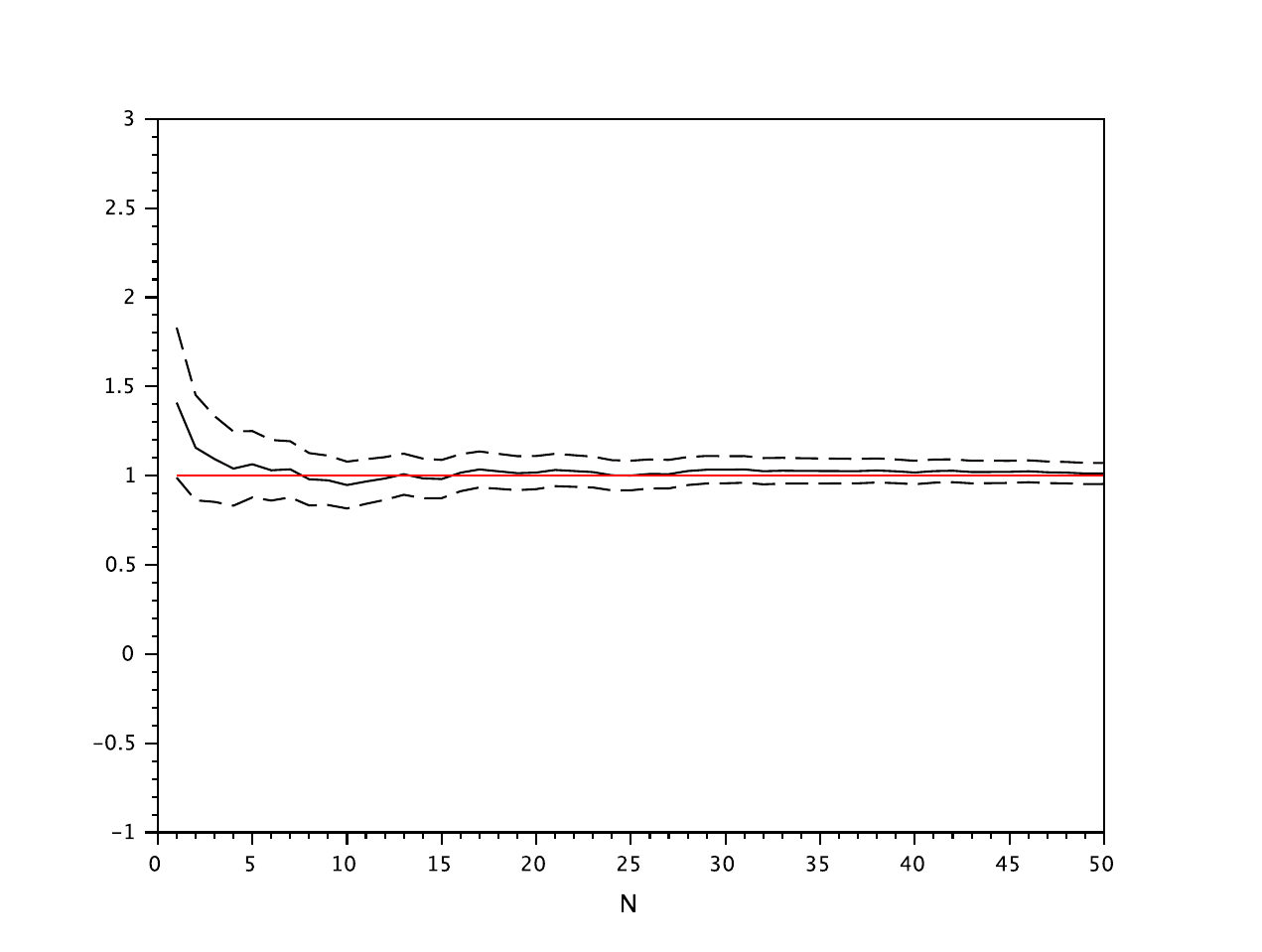} 
\end{tabular}
\end{center}
\caption{Plots of $N\mapsto\overline\theta_N$ (black line) and of the bounds of the $89\%$-ACIs (dashed black lines) for Model (B) with $H = 0.7$ (left) and $H = 0.9$ (right).}
\label{plots_Model_B}
\end{figure}
\begin{figure}[h!]
\begin{center}
\begin{tabular}{cc}
\includegraphics[scale=0.25]{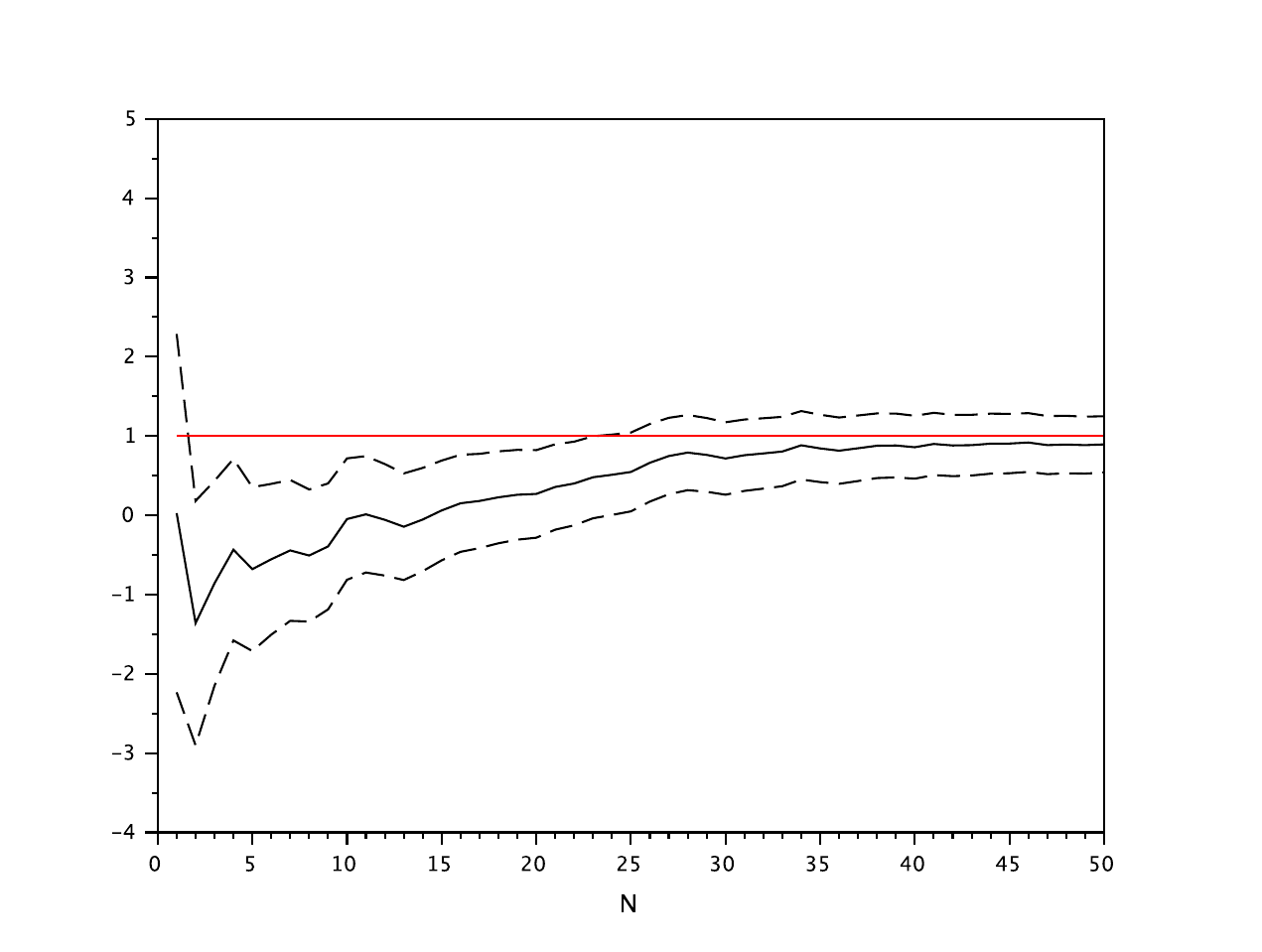} &
\includegraphics[scale=0.25]{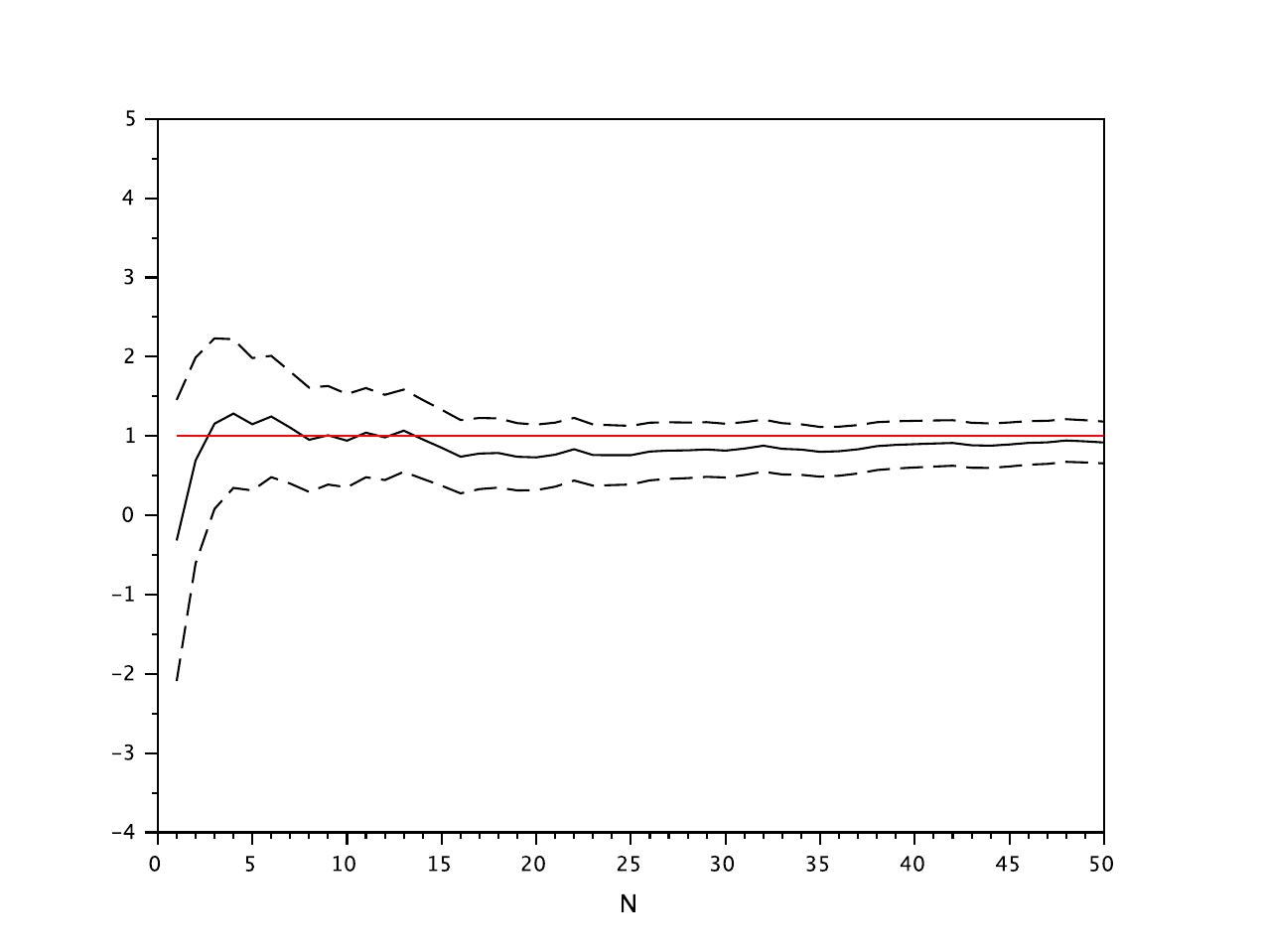} 
\end{tabular}
\end{center}
\caption{Plots of $N\mapsto\overline\theta_N$ (black line) and of the bounds of the $89\%$-ACIs (dashed black lines) for Model (C) with $H = 0.7$ (left) and $H = 0.9$ (right).}
\label{plots_Model_C}
\end{figure}
\newline
In Figures \ref{plots_Model_A}, \ref{plots_Model_B} and \ref{plots_Model_C}, for $N = 1,\dots,50$, $\overline\theta_N$ and the bounds of the $89\%$-asymptotic confidence interval (ACI) in Proposition \ref{ACI_approximation_LS_Young} (with $\alpha = 5\%$ and $\lambda = 0.11/(2\alpha) > 1$) are plotted for one of the 100 datasets generated from Models (A), (B) and (C) respectively. These figures illustrate both that $\overline\theta_N$ is consistent and its rate of convergence.
%


%
\section{Concluding remarks and perspectives}\label{conclusion_section}
This section summarizes the main insights gained throughout the paper and outlines potential directions for future research. As previously discussed, we investigated the estimation of the drift parameter based on $N$ independent copies of the solution to a stochastic differential equation driven by a multiplicative fractional Brownian noise with Hurst parameter $H\in (1/3,1)$. Our primary contributions encompass a rigorous analysis of the so-called fake-estimator, alongside a detailed study of the corrections required to transition from this theoretical construction to computable estimators. This transition involves replacing the Skorokhod integral with a suitably corrected pathwise integral.

For $H > 1/2$, this is achieved by reformulating the Malliavin derivative into a form that does not depend on $B$. Conversely, for $H\in (1/3,1/2]$, we establish a new result that allows us to handle the necessary corrections, yielding an expression analogous to the $H > 1/2$ regime - a result we believe is of independent mathematical interest. Leveraging these tools, we propose a computable fixed-point estimator and prove that it is well-defined. Furthermore, within the considered multiplicative framework, we establish both an asymptotic confidence interval and a non-asymptotic risk bound valid for any $H\in (1/3,1]$.\\

We now turn to possible extensions of our results, which open the way for future investigations.

A natural first consideration pertains to the parameter space. Our results can be naturally extended to the estimation of $\gamma_0\in\mathbb R$ when $\theta_0 =\texttt f(\gamma_0)$ and $\texttt f$ is a strictly monotonic continuous function from $\mathbb R$ into $(0,\infty)$. To that purpose, let us define the transformed estimator $\widehat\gamma_N =\texttt f^{-1}(\overline\theta_{N}^{\mathfrak c,\mathfrak d})$ of $\gamma_0$. On the one hand, given a confidence interval $\widehat I_N$ for $\theta_0$ derived from Proposition \ref{ACI_approximation_LS_Young} or \ref{ACI_approximation_LS_rough}, the random interval $\texttt f^{-1}(\widehat I_N)$ achieves the exact same coverage probability for $\gamma_0$. On the other hand, if $\texttt f^{-1}$ is Lipschitz continuous with constant $\texttt L > 0$, then the non-asymptotic risk bounds established in Propositions \ref{risk_bound_approximation_LS_Young} and \ref{risk_bound_approximation_LS_rough} extend to the transformed parameter as
\begin{displaymath}
\mathbb E(|\widehat\gamma_N -\gamma_0|^2)
\leqslant
\texttt L^2
\mathbb E(|\overline\theta_{N}^{\mathfrak c,\mathfrak d} -\theta_0|^2).
\end{displaymath}
This inequality guarantees that the original rate of convergence is fully preserved under the transformation.

Another compelling extension involves adapting our methodology to Equation (\ref{main_equation}) when $B$ is a $d$-dimensional fractional Brownian motion ($d\in\mathbb N^*$). To that purpose, consider $N$ independent copies
\begin{displaymath}
(X^1(x_0 +\varepsilon e_{\ell}))_{\ell = 0,\dots,d},
\dots,(X^N(x_0 +\varepsilon e_{\ell}))_{\ell = 0,\dots,d}
\quad {\rm of}\quad
(X(x_0 +\varepsilon e_{\ell}))_{\ell = 0,\dots,d},
\end{displaymath}
where $e_0 = 0$, $(e_1,\dots,e_d)$ is the canonical basis of $\mathbb R^d$, $\varepsilon > 0$ (fixed) and $X(\xi)$ is the solution of Equation (\ref{main_equation}) starting from $\xi\in\mathbb R^d$. First, note that a $d$-dimensional version of our fake-estimator is given by
\begin{displaymath}
\widehat\theta_N =
\left(\sum_{i = 1}^{N}\int_{0}^{T}b(X_{s}^{i}(x_0))^{\top}b(X_{s}^{i}(x_0))ds\right)^{-1}
\left(\sum_{i = 1}^{N}\int_{0}^{T}b(X_{s}^{i}(x_0))^{\top}\delta X_{s}^{i}(x_0)\right).
\end{displaymath}
Moreover, for any $i\in\{1,\dots,N\}$, let $\mathbf J_{u}^{i}$ ($u\in [0,T]$) be the Jacobian matrix of the map $\xi\mapsto X_{u}^{i}(\xi)$, and note that in the spirit of \cite{CM21}, Section 3.3, an estimator of $\mathbf J_{u}^{i}$ is given by
\begin{displaymath}
\widehat{\bf J}_{N,u}^{i} :=
\left(\frac{X_{u}^{i,j}(x_0 +\varepsilon e_{\ell}) - X_{u}^{i,j}}{\varepsilon}\right)_{j,\ell = 1,\dots,d}.
\end{displaymath}
Thus, since the relationship between the Skorokhod and pathwise integrals in Theorem 3.1 of \cite{ST22} holds for $d > 1$, one could investigate the computable estimator derived from $\widehat\theta_N$ by replacing $\mathbf J_{\cdot}^{i}$ by $\widehat{\bf J}_{N,.}^{i}$ in the following well-known equality:
\begin{displaymath}
\mathbf D_{s}^{i}X_{t}^{i} =
\mathbf 1_{[0,t]}(s)\mathbf J_{t}^{i}(\mathbf J_{s}^{i})^{-1}\sigma(X_{s}^{i}),
\quad s,t\in [0,T].
\end{displaymath}
%


%
\section{Proofs}\label{section_proofs}
This section is dedicated to the proofs of our main results. We begin by proving the probabilistic statements of Section \ref{section_probabilistic_preliminaries}. Subsequently, we proceed to the proofs of the statistical results. The proofs for both Hurst regimes follow a specific logical sequence. Recall that our approach builds upon an infeasible estimator $\widehat\theta_N$, which motivates the introduction of the computable estimator $\overline\theta_N = I_N + R_N$, where $R_N$ is defined as the fixed-point of a specific mapping $\Theta_N$.
\\
\\
{\bf The moderately irregular regime $(H > 1/2)$.} The proofs for this regime are structured as follows:
\begin{enumerate}
 \item {\bf Existence and uniqueness.} We first prove that, under the required conditions, the mapping $\Theta_N$ admits a unique fixed-point $R_N$. This theoretically justifies the introduction of the truncated estimator $\overline\theta_{N}^{\mathfrak c}$, which is restricted to the "good event" where this condition is automatically satisfied.
 \item {\bf Control of the error event.} We establish that the probability of the complementary "bad event" (where the fixed-point might fail to exist) vanishes asymptotically.
 \item {\bf Asymptotic normality.} Leveraging the control of the error event, we prove the asymptotic normality of $\overline\theta_{N}^{\mathfrak c}$, which allows us to construct valid asymptotic confidence intervals.
 \item {\bf Convergence rates.} Finally, in order to establish non-asymptotic convergence rates, we must strictly lower-bound the denominator of our estimator. We therefore prove the convergence rate for the doubly truncated estimator $\overline\theta_{N}^{\mathfrak c,\mathfrak d}$.
\end{enumerate}
{\bf The rough regime $(H < 1/2)$.} We then move to the rough framework, where the notation for the estimators remains identical. Despite the different stochastic calculus tools required, the logical roadmap mirrors the previous case:
\begin{enumerate}
 \item We prove the existence and uniqueness of the fixed-point under some conditions related to the rough paths framework, naturally leading to the same truncated estimator $\overline\theta_{N}^{\mathfrak c}$.
 \item We construct the asymptotic confidence interval for $\overline\theta_{N}^{\mathfrak c}$ by showing that the error event probability vanishes appropriately.
 \item Finally, we omit the formal proof of the non-asymptotic convergence rate for the doubly truncated estimator $\overline\theta_{N}^{\mathfrak c,\mathfrak d}$. As previously noted, the underlying mathematical arguments are identical to those detailed for the $H > 1/2$ regime.
\end{enumerate}
%


%
\subsection{Proof of Proposition \ref{Malliavin_derivative_X}}\label{section_proof_Malliavin_derivative_X}
It is well-known (see, for instance, Theorem 2.2.1 in \cite{NUALART06}) that for any $s,t\in [0,T]$ with $s < t$,
\begin{displaymath}
\mathbf D_sX_t =
\sigma(X_s) +
\theta_0\int_{s}^{t}b'(X_u)\mathbf D_sX_udu +
\int_{s}^{t}\sigma'(X_u)\mathbf D_sX_udB_u.
\end{displaymath}
Then, by the change-of-variable formula for the Young integral (resp. rough integral) when $H\in (1/2,1)$ (resp. $H\in (1/3,1/2]$),
\begin{eqnarray*}
 \mathbf D_sX_t & = &
 \sigma(X_s)\exp\left(
 \theta_0\int_{s}^{t}b'(X_u)du +
 \int_{s}^{t}\sigma'(X_u)dB_u\right)\\
 & = &
 \sigma(X_s)\exp\left(
 \theta_0\int_{s}^{t}b'(X_u)du +
 \int_{s}^{t}\frac{\sigma'(X_u)}{\sigma(X_u)}
 (\theta_0b(X_u)du +\sigma(X_u)dB_u -\theta_0b(X_u)du)\right)\\
 & = &
 \sigma(X_s)\exp\left(
 \theta_0\int_{s}^{t}\left(b'(X_u) -\frac{\sigma'(X_u)b(X_u)}{\sigma(X_u)}\right)du +
 \int_{s}^{t}\frac{\sigma'(X_u)}{\sigma(X_u)}dX_u\right).
\end{eqnarray*}
Moreover, once again by using the change-of-variable formula,
\begin{displaymath}
\int_{s}^{t}\frac{\sigma'(X_u)}{\sigma(X_u)}dX_u =
\log(|\sigma(X_t)|) -\log(|\sigma(X_s)|).
\end{displaymath}
Therefore,
\begin{eqnarray*}
 \mathbf D_sX_t
 & = &
 \sigma(X_s)\exp(\log(|\sigma(X_t)|) -\log(|\sigma(X_s)|))
 \exp\left(
 \theta_0\int_{s}^{t}\left(b'(X_u) -
 \frac{\sigma'(X_u)b(X_u)}{\sigma(X_u)}\right)du\right)\\
 & = &
 \sigma(X_t)\exp\left(
 \theta_0\int_{s}^{t}\left(b'(X_u) -
 \frac{\sigma'(X_u)b(X_u)}{\sigma(X_u)}\right)du\right)
 \hspace{0.75cm}\Box
\end{eqnarray*}


%
\subsection{Proof of Proposition \ref{2D_Young_integral_R}}
First of all, let us define the rectangular increments of $R$:
\begin{displaymath}
\Delta_{(s,t),(u,v)}R := R(s,t) - R(s,v) - (R(u,t) - R(u,v)),
\quad s,t,u,v\in [0,T].
\end{displaymath}
Moreover, let $\pi =\pi_n = (t_0,\dots,t_n)$ be a dissection of $[0,T]$ such that
\begin{displaymath}
|\pi_n| :=
\max_{j\in\{0,\dots,n - 1\}}|t_{j + 1} - t_j|\rightarrow 0
\quad {\rm when}\quad n\rightarrow\infty.
\end{displaymath}
By \cite{FV10}, Section 6.4, the Riemann sum
\begin{displaymath}
I^{\pi_n}(x,R) =
\sum_{j = 0}^{n - 1}\sum_{i = 0}^{j}
x(t_i,t_j)\Delta_{(t_i,t_j),(t_{i + 1},t_{j + 1})}R
\end{displaymath}
converges to the 2D Young integral of $x\mathbf 1_{\Delta_T}$ with respect to $R$ when $n\rightarrow\infty$. So, it is sufficient to establish that
\begin{displaymath}
\lim_{n\rightarrow\infty}I^{\pi_n}(x,R) =
\alpha_H\int_{0}^{T}\int_{0}^{t}x(s,t)|t - s|^{2H - 2}dsdt
\quad\textrm{to conclude}.
\end{displaymath}
To that purpose, note that
\begin{displaymath}
I^{\pi}(x,R) = I_{<}^{\pi}(x,R) + I_{=}^{\pi}(x,R),
\end{displaymath}
where
\begin{eqnarray*}
 I_{<}^{\pi}(x,R) & = &
 \sum_{j = 0}^{n - 1}\sum_{i = 0}^{j - 2}x(t_i,t_j)\Delta_{(t_i,t_j),(t_{i + 1},t_{j + 1})}R\\
 & &
 \hspace{2cm}{\rm and}\quad
 I_{=}^{\pi}(x,R) =
 \sum_{j = 0}^{n - 1}\sum_{i = j - 1}^{j}x(t_i,t_j)\Delta_{(t_i,t_j),(t_{i + 1},t_{j + 1})}R.
\end{eqnarray*}
The proof of Lemma \ref{2D_Young_integral_R} is dissected in three steps.
\\
\\
{\bf Step 1 (preliminaries).} Let $\partial_2R$ be the partial derivative of $R$ with respect to its second variable. First, note that for every $i,j\in\{1,\dots,n - 1\}$ such that $i\neq j$,
\begin{equation}\label{2D_Young_integral_R_3}
\Delta_{(t_i,t_j),(t_{i + 1},t_{j + 1})}R =
\int_{t_j}^{t_{j + 1}}(\partial_2R(t_{i + 1},t) -\partial_2R(t_i,t))dt
\end{equation}
and
\begin{equation}\label{2D_Young_integral_R_4}
\Delta_{(t_i,t_j),(t_{i + 1},t_{j + 1})}R =
\frac{1}{2}(|t_{i + 1} - t_j|^{2H} + |t_{j + 1} - t_i|^{2H} - |t_{i + 1} - t_{j + 1}|^{2H} - |t_i - t_j|^{2H}).
\end{equation}
Now, for every $(s,t)\in\Delta_T$,
\begin{eqnarray}
 \partial_2R(s,t) & = &
 H(t^{2H - 1} - (t - s)^{2H - 1})
 \label{2D_Young_integral_R_5}
 =
 H(2H - 1)\int_{0}^{s}|t - u|^{2H - 2}du.
\end{eqnarray}
Finally, since $z\mapsto z^{2H}$ is a $2H$-H\"older continuous function from $\mathbb R_+$ into itself,
\begin{equation}\label{2D_Young_integral_R_6}
|t - s|^{2H}\leqslant
||t - s|^{2H} - |r - s|^{2H}| + |r - s|^{2H}\leqslant
|t - r|^{2H} + |s - r|^{2H}.
\end{equation}
{\bf Step 2 (control of $I_{=}^{\pi}(x,R)$).} By the condition (\ref{2D_Young_integral_R_1}),
\begin{displaymath}
|I_{=}^{\pi}(x,R)|\leqslant
\mathfrak c_x\sum_{j = 1}^{n - 1}|t_j - t_{j - 1}|^{\alpha}|\Delta_{(t_{j - 1},t_j),(t_j,t_{j + 1})}R|.
\end{displaymath}
Moreover, by Equality (\ref{2D_Young_integral_R_4}) and Inequality (\ref{2D_Young_integral_R_6}),
\begin{eqnarray*}
 |\Delta_{(t_{j - 1},t_j),(t_j,t_{j + 1})}R| & = &
 \frac{1}{2}||t_{j + 1} - t_j + t_j - t_{j - 1}|^{2H} - |t_j - t_{j + 1}|^{2H} - |t_{j - 1} - t_j|^{2H}|\\
 & \leqslant &
 |t_{j + 1} - t_j|^{2H} + |t_j - t_{j - 1}|^{2H}.
\end{eqnarray*}
So, by the H\"older inequality with the conjugate exponents $p = (1 +\varepsilon)/\alpha$ - where $\varepsilon = \alpha + 2H - 1$ - and $q = 1 + 1/(p - 1)$,
\begin{eqnarray*}
 |I_{=}^{\pi}(x,R)| & \leqslant &
 \mathfrak c_x\sum_{j = 1}^{n - 1}|t_j - t_{j - 1}|^{\alpha}(|t_{j + 1} - t_j|^{2H} + |t_j - t_{j - 1}|^{2H})\\
 & \leqslant &
 \mathfrak c_x
 T\max_{j\in\{1,\dots,n - 1\}}|t_j - t_{j - 1}|^{\alpha + 2H - 1}\\
 & &
 \hspace{1.5cm} +
 \mathfrak c_x
 \left(T\max_{j\in\{1,\dots,n - 1\}}|t_{j + 1} - t_j|^{\varepsilon}\right)^{\frac{1}{p}}
 \left(T\max_{j\in\{1,\dots,n - 1\}}|t_j - t_{j - 1}|^{2Hq - 1}\right)^{\frac{1}{q}}.
\end{eqnarray*}
Therefore, since $2Hq - 1 =\alpha + 2H - 1 =\varepsilon$ and $1/p + 1/q = 1$,
\begin{displaymath}
|I_{=}^{\pi}(x,R)|\leqslant
2\mathfrak c_xT
\max_{j\in\{1,\dots,n\}}|t_j - t_{j - 1}|^{\alpha + 2H - 1}
\xrightarrow[n\rightarrow\infty]{} 0.
\end{displaymath}
{\bf Step 3 (control of $I_{<}^{\pi}(x,R)$).} By Equalities (\ref{2D_Young_integral_R_3}) and (\ref{2D_Young_integral_R_5}),
\begin{eqnarray*}
 I_{<}^{\pi}(x,R) & = &
 \sum_{j = 0}^{n - 1}\sum_{i = 0}^{j - 2}
 x(t_i,t_j)\int_{t_j}^{t_{j + 1}}(\partial_2R(t_{i + 1},t) -\partial_2R(t_i,t))dt\\
 & = &
 \alpha_H\sum_{j = 0}^{n - 1}\sum_{i = 0}^{j - 2}
 x(t_i,t_j)\int_{t_j}^{t_{j + 1}}\int_{t_i}^{t_{i + 1}}|t - s|^{2H - 2}dsdt = U_n + V_n
\end{eqnarray*}
with
\begin{eqnarray*}
 U_n & = &
 \alpha_H\sum_{j = 0}^{n - 1}\sum_{i = 0}^{j - 2}
 \int_{t_j}^{t_{j + 1}}\int_{t_i}^{t_{i + 1}}x(s,t)|t - s|^{2H - 2}dsdt\\
 & &
 \hspace{2.5cm}{\rm and}\quad
 V_n =
 \alpha_H\sum_{j = 0}^{n - 1}\sum_{i = 0}^{j - 2}
 \int_{t_j}^{t_{j + 1}}\int_{t_i}^{t_{i + 1}}(x(t_i,t_j) - x(s,t))|t - s|^{2H - 2}dsdt.
\end{eqnarray*}
On the one hand, by the condition (\ref{2D_Young_integral_R_1}),
\begin{eqnarray*}
 & &
 \left|U_n -\alpha_H\int_{0}^{T}\int_{0}^{t}x(s,t)|t - s|^{2H - 2}dsdt\right|\\
 & &
 \hspace{2cm} =
 \left|\alpha_H\sum_{j = 0}^{n - 1}
 \int_{t_j}^{t_{j + 1}}\int_{t_{j - 1}}^{t}x(s,t)|t - s|^{2H - 2}dsdt\right|\\
 & &
 \hspace{2cm}\leqslant
 \frac{\mathfrak c_x|\alpha_H|}{\alpha + 2H - 1}\sum_{j = 0}^{n - 1}
 \int_{t_j}^{t_{j + 1}}(t - t_{j - 1})^{\alpha + 2H - 1}dt\\
 & &
 \hspace{2cm}\leqslant
 \frac{\mathfrak c_x|\alpha_H|}{(\alpha + 2H - 1)(\alpha + 2H)}\sum_{j = 0}^{n - 1}
 (t_{j + 1} - t_{j - 1})^{\alpha + 2H}\lesssim
 |\pi_n|^{\alpha + 2H - 1}
 \xrightarrow[n\rightarrow\infty]{} 0.
\end{eqnarray*}
On the other hand, by the condition (\ref{2D_Young_integral_R_2}),
\begin{eqnarray*}
 |V_n| & \leqslant &
 \mathfrak c_x|\alpha_H|\sum_{j = 0}^{n - 1}\sum_{i = 0}^{j - 2}
 \int_{t_j}^{t_{j + 1}}\int_{t_i}^{t_{i + 1}}(|t_i - s|^{\alpha} + |t_j - t|^{\alpha})|t - s|^{2H - 2}dsdt\\
 & \leqslant &
 \mathfrak c_x|\alpha_H|\sum_{i = 0}^{n - 3}
 \int_{t_i}^{t_{i + 1}}|t_i - s|^{\alpha}\int_{t_{i + 1}}^{T}(t - s)^{2H - 2}dtds\\
 & &
 \hspace{2cm} +
 \mathfrak c_x|\alpha_H|\sum_{j = 0}^{n - 1}
 \int_{t_j}^{t_{j + 1}}|t_j - t|^{\alpha}\int_{0}^{t_{j - 1}}(t - s)^{2H - 2}dsdt\\
 & \leqslant &
 \mathfrak c_x\sum_{i = 0}^{n - 3}
 \int_{t_i}^{t_{i + 1}}(s - t_i)^{\alpha}(t_{i + 1} - s)^{2H - 1}ds\\
 & &
 \hspace{2cm} +
 \mathfrak c_x\sum_{j = 0}^{n - 1}
 \int_{t_j}^{t_{j + 1}}(t - t_j)^{\alpha}(t - t_{j - 1})^{2H - 1}dt
 \leqslant\mathfrak c_x|\pi_n|^{\alpha + 2H - 1}\xrightarrow[n\rightarrow\infty]{} 0.
\end{eqnarray*}
Therefore,
\begin{displaymath}
\lim_{n\rightarrow\infty}I_{<}^{\pi}(x,R) =
\alpha_H\int_{0}^{T}\int_{0}^{t}x(s,t)|t - s|^{2H - 2}dsdt
\hspace{0.75cm}\Box
\end{displaymath}
%


%
\subsection{Proof of Proposition \ref{relationship_Skorokhod_rough}}\label{section_proof_relationship_Skorokhod_rough}
In the sequel, the set of all the dissections of $[0,T]$ is denoted by $\mathfrak D_T$. First of all, let us define the 1D and 2D $p$-variation norms ($p\in [1,\infty)$). A continuous function $x : [0,T]\rightarrow\mathbb R$ is of finite $p$-variation if and only if
\begin{displaymath}
\|x\|_{p\textrm{-var},T} :=
\left(\sup_{(t_j)\in\mathfrak D_T}
\sum_j|x(t_{j + 1}) - x(t_j)|^p\right)^{\frac{1}{p}} <\infty,
\end{displaymath}
and a continuous function $\rho : [0,T]^2\rightarrow\mathbb R$ is of finite $p$-variation if and only if
\begin{displaymath}
\|\rho\|_{p\textrm{-var},[0,T]^2} :=
\left(\sup_{(s_i),(t_j)\in\mathfrak D_T}
\sum_{i,j}|\rho(t_{j + 1},s_i) -\rho(t_j,s_i) -
(\rho(t_{j + 1},s_{i + 1}) -\rho(t_j,s_{i + 1}))|^p\right)^{\frac{1}{p}} <\infty.
\end{displaymath}
The proof of Proposition \ref{relationship_Skorokhod_rough} relies on Proposition \ref{2D_Young_integral_R} and on the following technical lemma.
%


%
\begin{lemma}\label{properties_L}
For every $s,t\in [0,T]$, consider
\begin{displaymath}
\overline L(s,t) =
\varphi(X_t)\exp\left(\theta_0\int_{s}^{t}\psi(X_u)du\right)\mathbf 1_{[0,t)}(s)
\quad\textrm{and}\quad
L(s,t) =
\overline L(s,t) -\varphi(X_t)\mathbf 1_{[0,t)}(s).
\end{displaymath}
Under the assumptions of Proposition \ref{relationship_Skorokhod_rough},
\begin{enumerate}
 \item There exists a deterministic constant $\mathfrak c_{\ref{properties_L}} > 0$ such that
 \begin{displaymath}
 |L(s,t)|\leqslant
 \mathfrak c_{\ref{properties_L}}|t - s|,
 \quad\forall (s,t)\in\Delta_T,
 \end{displaymath}
 and for every $\alpha\in (0,H)$,
 \begin{displaymath}
 |L(s,t) - L(u,v)|\leqslant
 \mathfrak c_{\ref{properties_L}}(1\vee\|X\|_{\alpha,T})(|t - v|^{\alpha} + |s - u|),
 \quad\forall (s,t),(u,v)\in\Delta_T,
 \end{displaymath}
 where $\|X\|_{\alpha,T}$ is the $\alpha$-H\"older norm of $X$ over $[0,T]$.
 \item For every $p > 1/H$, the random variables $\|\overline L(0,\cdot)\|_{p\normalfont{\textrm{-var}},T}$ and $\|\overline L\|_{p\normalfont{\textrm{-var}},[0,T]^2}$ belong to $\mathbb L^2(\Omega)$.
\end{enumerate}
\end{lemma}
The proof of Lemma \ref{properties_L} is postponed to Section \ref{section_proof_properties_L}. Note that Lemma \ref{properties_L}.(1) allows to apply our Proposition \ref{2D_Young_integral_R} to the $L(\omega)$'s, and that Lemma \ref{properties_L}.(2) allows to apply \cite{ST22}, Theorem 3.1 to $Y =\pi\circ X$.
\\
\\
First, as established in \cite{FH14}, Section 7.3, the paths of $Y$ are controlled by those of $B$. Precisely, for every $(s,t)\in\Delta_T$,
\begin{displaymath}
Y_t - Y_s = Y_s'(B_t - B_s) +\mathfrak R(s,t),
\end{displaymath}
where
\begin{displaymath}
Y_s' =\pi'(X_s)X_s' =\pi'(X_s)\sigma(X_s) =\varphi(X_s)
\end{displaymath}
is the Gubinelli derivative of $Y$ at time $s$, and $\mathfrak R$ is a stochastic process which paths are continuous and of finite $p/2$-variation for every $p\in (1/H,3)$. Moreover, recalling that $\varphi$ and $\psi$ are defined in (\ref{definition_phi_psi}), by Proposition \ref{Malliavin_derivative_X}, for every $(s,t)\in\Delta_T$,
\begin{displaymath}
\mathbf D_sY_t =
\overline L(s,t) = L(s,t) + Y_t'.
\end{displaymath}
Let us introduce $V(s) := R(s,s)$ ($s\in [0,T]$). Thanks to Lemma \ref{properties_L}.(2), by \cite{ST22}, Theorem 3.1,
\begin{eqnarray*}
 \int_{0}^{T}Y_s\delta B_s
 & = &
 \int_{0}^{T}Y_sdB_s -\frac{1}{2}\int_{0}^{T}Y_s'dV(s) -
 \int_{0 < s < t < T}
 (\mathbf D_sY_t - Y_t')dR(s,t)\\
 & = &
 \int_{0}^{T}Y_sdB_s - H\int_{0}^{T}\varphi(X_s)s^{2H - 1}ds -
 \int_{0 < s < t < T}L(s,t)dR(s,t).
\end{eqnarray*}
Now, thanks to Lemma \ref{properties_L}.(1), by Proposition \ref{2D_Young_integral_R},
\begin{eqnarray*}
 \int_{0 < s < t < T}L(s,t)dR(s,t)
 & = &
 \alpha_H\int_{0}^{T}\int_{0}^{t}
 L(s,t)|t - s|^{2H - 2}dsdt\\
 & = &
 \alpha_H\int_{0}^{T}\int_{0}^{t}
 \varphi(X_t)\left(\exp\left(\theta_0\int_{s}^{t}\psi(X_u)du\right) - 1\right)|t - s|^{2H - 2}dsdt.
\end{eqnarray*}
Therefore,
\begin{eqnarray*}
 \int_{0}^{T}Y_s\delta B_s
 & = &
 \int_{0}^{T}Y_sdB_s - H\int_{0}^{T}\varphi(X_s)s^{2H - 1}ds\\
 & &
 \hspace{2cm} -\alpha_H\int_{0}^{T}\int_{0}^{t}
 \varphi(X_t)\left(\exp\left(\theta_0\int_{s}^{t}\psi(X_u)du\right) - 1\right)|t - s|^{2H - 2}dsdt,
\end{eqnarray*}
leading to
\begin{eqnarray*}
 \int_{0}^{T}b(X_s)\delta X_s
 & = &
 \int_{0}^{T}b(X_s)dX_s - H\int_{0}^{T}\varphi(X_s)s^{2H - 1}ds\\
 & &
 \hspace{2cm} -\alpha_H\int_{0}^{T}\int_{0}^{t}
 \varphi(X_t)\left(\exp\left(\theta_0\int_{s}^{t}\psi(X_u)du\right) - 1\right)|t - s|^{2H - 2}dsdt
 \hspace{0.75cm}\Box
\end{eqnarray*}
%


%
\subsubsection{Proof of Lemma \ref{properties_L}}\label{section_proof_properties_L}
For every $(s,t)\in\Delta_T$, consider
\begin{displaymath}
\overline\lambda(s,t) :=
\exp\left(\theta_0\int_{s}^{t}\psi(X_u)du\right)
\quad {\rm and}\quad
\lambda(s,t) :=
\overline\lambda(s,t) - 1.
\end{displaymath}
\begin{enumerate}
 \item First, for every $(s,t)\in\Delta_T$,
 \begin{eqnarray*}
  |L(s,t)| & = &
  |\varphi(X_t)\lambda(s,t)|\\
  & \leqslant &
  \mathfrak c_1|t - s|
  \quad {\rm with}\quad
  \mathfrak c_1 =
  |\theta_0|\cdot
  \|\varphi\|_{\infty}\|\psi\|_{\infty}e^{|\theta_0|\cdot\|\psi\|_{\infty}T}.
 \end{eqnarray*}
 Now, for any $\alpha\in (0,H)$ and $(s,t),(u,v)\in\Delta_T$,
 \begin{eqnarray*}
  |L(s,t) - L(u,v)| & = &
  |\varphi(X_t)\lambda(s,t) -\varphi(X_v)\lambda(u,v)|\\
  & \leqslant &
  |\varphi(X_t)|\cdot |\lambda(s,t) -\lambda(u,v)| +
  |\lambda(u,v)|\cdot |\varphi(X_t) -\varphi(X_v)|\\
  & \leqslant &
  \|\varphi\|_{\infty}|\lambda(s,t) -\lambda(u,v)| +
  \|\varphi'\|_{\infty}\|X\|_{\alpha,T}
  (1 + e^{|\theta_0|\cdot\|\psi\|_{\infty}T})|t - v|^{\alpha}.
 \end{eqnarray*}
 Moreover,
 \begin{eqnarray*}
  |\lambda(s,t) -\lambda(u,v)| & \leqslant &
  \mathfrak c_2
  \left|\int_{s}^{t}\psi(X_r)dr -\int_{u}^{v}\psi(X_r)dr\right|
  \quad {\rm with}\quad
  \mathfrak c_2 =
  |\theta_0|e^{|\theta_0|\cdot\|\psi\|_{\infty}T}\\
  & = &
  \mathfrak c_2
  \left|\int_{s}^{u}\psi(X_r)dr -\int_{t}^{v}\psi(X_r)dr\right|
  \leqslant
  \mathfrak c_2\|\psi\|_{\infty}(|s - u| + |t - v|).
 \end{eqnarray*}
 Then,
 \begin{equation}\label{properties_L_1}
 |L(s,t) - L(u,v)|\leqslant
 \mathfrak c_3(1\vee\|X\|_{\alpha,T})(|t - v|^{\alpha} + |s - u|)
 \end{equation}
 with
 \begin{displaymath}
 \mathfrak c_3 =
 (\mathfrak c_2\|\varphi\|_{\infty}\|\psi\|_{\infty})\vee
 [\|\varphi'\|_{\infty}(1 + e^{|\theta_0|\cdot\|\psi\|_{\infty}T}) +
 \mathfrak c_2\|\varphi\|_{\infty}\|\psi\|_{\infty}T^{1 -\alpha}].
 \end{displaymath}
 \item Consider $p > 1/H$. First, by Inequality (\ref{properties_L_1}), for every $s,t\in [0,T]$,
 \begin{eqnarray*}
  |\overline L(0,t) -\overline L(0,s)| & \leqslant &
  |L(0,t) - L(0,s)| + |\varphi(X_t) -\varphi(X_s)|\\
  & \leqslant &
  \mathfrak c_4(1\vee\|X\|_{\alpha,T})|t - s|^{\alpha}
  \quad {\rm with}\quad
  \mathfrak c_4 =\mathfrak c_3 +\|\varphi'\|_{\infty}
  \quad {\rm and}\quad
  \alpha =\frac{1}{p}.
 \end{eqnarray*}
 Then,
 \begin{eqnarray*}
  \|\overline L(0,\cdot)\|_{p\textrm{-var},T}^{p} & = &
  \sup_{(t_j)\in\mathfrak D_T}\sum_j|\overline L(0,t_{j + 1}) -\overline L(0,t_j)|^p\\
  & \leqslant &
  \mathfrak c_{4}^{p}(1\vee\|X\|_{\alpha,T})^p
  \sup_{(t_j)\in\mathfrak D_T}\sum_j|t_{j + 1} - t_j|
  \leqslant
  \mathfrak c_{4}^{p}T(1\vee\|X\|_{\alpha,T})^p.
 \end{eqnarray*}
 Now, using the definitions of $\overline L$, $\overline\lambda$ and $\lambda$, for every $s,t,u,v\in [0,T]$ such that $s,u < t,v$,
 \begin{eqnarray*}
  & &
  \overline L(s,t) -\overline L(u,t) - (\overline L(s,v) -\overline L(u,v))\\
  & &
  \hspace{2cm} =
  \varphi(X_t)(\overline\lambda(s,t) -\overline\lambda(u,t)) -
  \varphi(X_v)(\overline\lambda(s,v) -\overline\lambda(u,v))\\
  & &
  \hspace{2cm} =
  \varphi(X_t)(\overline\lambda(s,u)\overline\lambda(u,t) -\overline\lambda(u,t)) -
  \varphi(X_v)(\overline\lambda(s,u)\overline\lambda(u,v) -\overline\lambda(u,v))\\
  & &
  \hspace{2cm} =
  (\overline\lambda(s,u) - 1)(\varphi(X_t)\overline\lambda(u,t) -\varphi(X_v)\overline\lambda(u,v))\\
  & &
  \hspace{2cm} =
  \lambda(s,u)(\varphi(X_t)(\overline\lambda(u,t) -\overline\lambda(u,v)) +
  \overline\lambda(u,v)(\varphi(X_t) -\varphi(X_v)))\\
  & &
  \hspace{2cm} =
  \lambda(s,u)\overline\lambda(u,v)
  (\varphi(X_t)\lambda(v,t) +\varphi(X_t) -\varphi(X_v)),
 \end{eqnarray*}
 leading to
 \begin{eqnarray*}
  & &
  |\overline L(s,t) -\overline L(u,t) - (\overline L(s,v) -\overline L(u,v))|\\
  & &
  \hspace{1cm}\leqslant
  \|\overline\lambda\|_{\infty}|\lambda(s,u)|
  (\|\varphi\|_{\infty}|\lambda(v,t)| +\|\varphi'\|_{\infty}\|X\|_{\alpha,T}|t - v|^{\alpha})\\
  & &
  \hspace{1cm}\leqslant
  \mathfrak c_5|s - u|
  (\mathfrak c_1|t - v| +
  \|\varphi'\|_{\infty}\|X\|_{\alpha,T}|t - v|^{\alpha})
  \quad {\rm with}\quad
  \mathfrak c_5 = |\theta_0|\cdot\|\psi\|_{\infty}e^{2|\theta_0|\cdot\|\psi\|_{\infty}T}\\
  & &
  \hspace{1cm}\leqslant
  \mathfrak c_6(1\vee\|X\|_{\alpha,T})|s - u|^{\alpha}|t - v|^{\alpha}\\
  & &
  \hspace{3cm}{\rm with}\quad
  \mathfrak c_6 =\mathfrak c_5T^{1 -\alpha}(\mathfrak c_1T^{1 -\alpha} +\|\varphi'\|_{\infty})
  \quad\textrm{and (again)}\quad
  \alpha = p^{-1}.
 \end{eqnarray*}
 Then,
 \begin{eqnarray*}
  \|\overline L\|_{p\textrm{-var},[0,T]^2}^{p} & = &
  \sup_{(s_i),(t_j)\in\mathfrak D_T}\sum_{i,j}
  |\overline L(t_{j + 1},s_i) -\overline L(t_j,s_i) -
  (\overline L(t_{j + 1},s_{i + 1}) -\overline L(t_j,s_{i + 1}))|^p\\
  & \leqslant &
  \mathfrak c_{6}^{p}(1\vee\|X\|_{\alpha,T})^p
  \sup_{(s_i),(t_j)\in\mathfrak D_T}
  \sum_{i,j}|s_{i + 1} - s_i|\cdot |t_{j + 1} - t_j|
  \leqslant
  \mathfrak c_{6}^{p}T^2(1\vee\|X\|_{\alpha,T})^p.
 \end{eqnarray*}
 In conclusion, since $b$ and $\sigma$ are bounded, by \cite{FH14}, Proposition 8.3, and by \cite{FV10}, Theorem 15.33.(iii),
 \begin{eqnarray*}
  \mathbb E(\|\overline L(0,\cdot)\|_{p\textrm{-var},T}^{2})
  & \leqslant &
  \mathfrak c_{4}^{2}T^{2\alpha}(1 +\mathbb E(\|X\|_{\alpha,T}^{2})) <\infty\\
  & &
  \hspace{0.75cm}{\rm and}\quad
  \mathbb E(\|\overline L\|_{p\textrm{-var},[0,T]^2}^{2})
  \leqslant
  \mathfrak c_{6}^{2}T^{4\alpha}(1 +\mathbb E(\|X\|_{\alpha,T}^{2})) <\infty
  \hspace{0.75cm}\Box
 \end{eqnarray*}
\end{enumerate}
%


%
\subsection{Proof of Proposition \ref{existence_uniqueness_R_N_Young}}\label{section_proof_existence_uniqueness_R_N_Young}
Since $\varphi\leqslant 0$, $\Theta_N$ is nonnegative and, in particular, $\Theta_N(\mathbb R_+)\subset\mathbb R_+$. Moreover, by (\ref{existence_uniqueness_R_N_Young_1}), for every $r,\overline r\in\mathbb R_+$,
\begin{eqnarray*}
 |\Theta_N(\overline r) -\Theta_N(r)| & \leqslant &
 \frac{\alpha_H}{NTD_N}\sum_{i = 1}^{N}
 \int_{0}^{T}\int_{0}^{t}|t - s|^{2H - 2}|\varphi(X_{t}^{i})|
 \exp\left(I_N\int_{s}^{t}\psi(X_{u}^{i})du\right)\\
 & &
 \hspace{4cm}\times
 \left|\exp\left(\overline r\int_{s}^{t}\psi(X_{u}^{i})du\right) -
 \exp\left(r\int_{s}^{t}\psi(X_{u}^{i})du\right)\right|dsdt\\
 & \leqslant &
 \frac{\alpha_H\|\varphi\|_{\infty}}{NTD_N}
 e^{\|\psi\|_{\infty}|I_N|T}
 \sum_{i = 1}^{N}\int_{0}^{T}\int_{0}^{t}|t - s|^{2H - 2}\sup_{x\in\mathbb R_-}e^x\left|
 (\overline r - r)\int_{s}^{t}\psi(X_{u}^{i})du\right|dsdt\\
 & \leqslant &
 \overline\alpha_H\|\varphi\|_{\infty}\|\psi\|_{\infty}
 T^{2H}\frac{M_N}{D_N}|\overline r - r|
 \leqslant\mathfrak c|\overline r - r|.
\end{eqnarray*}
So, $\Theta_N$ is a contraction from $\mathbb R_+$ into itself, and then $R_N$ exists and is unique by Picard's fixed-point theorem.$\hspace{0.75cm}\Box$
%


%
\subsection{Proof of Proposition \ref{deviation_bound}}\label{section_proof_deviation_bound}
Note that
\begin{eqnarray*}
 \mathbb P(\Delta_{N}^{c}) & = &
 \mathbb P\left(\frac{M_N}{D_N} >
 \frac{\mathfrak c}{\overline\alpha_HT^{2H}\|\varphi\|_{\infty}\|\psi\|_{\infty}},
 D_N\geqslant\mathfrak e\right)\\
 & &
 \hspace{2cm} +
 \mathbb P\left(\frac{M_N}{D_N} >
 \frac{\mathfrak c}{\overline\alpha_HT^{2H}\|\varphi\|_{\infty}\|\psi\|_{\infty}},
 D_N <\mathfrak e\right)
 \quad {\rm with}\quad\mathfrak e =\frac{\|b\|_{f}^{2}}{2}\\
 & \leqslant &
 \mathbb P(\Delta_{N}^{\mathfrak e}) +
 \mathbb P(D_N <\mathfrak e)
 \quad {\rm with}\quad
 \Delta_{N}^{\mathfrak e} =\left\{
 M_N >\frac{\mathfrak c\mathfrak e}{
 \overline\alpha_HT^{2H}\|\varphi\|_{\infty}\|\psi\|_{\infty}}\right\}\cap
 \{D_N\geqslant\mathfrak e\}.
\end{eqnarray*}
On the one hand, using Markov's inequality and that $X^1,\dots,X^N$ are i.i.d. processes, we obtain
\begin{eqnarray}
 \mathbb P(D_N <\mathfrak e)
 & \leqslant &
 \mathbb P\left(|D_N -\|b\|_{f}^{2}| >\frac{\|b\|_{f}^{2}}{2}\right)
 \nonumber\\
 \label{deviation_bound_2}
 & \leqslant &
 \frac{1}{\mathfrak e^2}\mathbb E(|D_N -\|b\|_{f}^{2}|^2) =
 \frac{1}{\mathfrak e^2}\cdot\frac{1}{N^2T^2}
 {\rm var}\left(\sum_{i = 1}^{N}\int_{0}^{T}b(X_{s}^{i})^2ds\right)\\
 & = &
 \frac{1}{\mathfrak e^2N}{\rm var}\left(\frac{1}{T}\int_{0}^{T}b(X_s)^2ds\right)
 \leqslant
 \frac{\|b^2\|_{f}^{2}}{\mathfrak e^2N}.
 \nonumber
\end{eqnarray}
On the other hand,
\begin{eqnarray*}
 \mathbb P(\Delta_{N}^{\mathfrak e}) & \leqslant &
 \mathbb P\left[
 \frac{1}{N}\left|\sum_{i = 1}^{N}(\texttt b(X_{T}^{i}) -\texttt b(x_0))\right| >
 \log\left(\frac{\mathfrak c\mathfrak e}{
 \overline\alpha_HT^{2H}\|\varphi\|_{\infty}\|\psi\|_{\infty}}\right)
 \frac{\mathfrak e}{\|\psi\|_{\infty}}\right]\\
 & \leqslant &
 \mathbb P\left[
 \frac{1}{N}\left|\sum_{i = 1}^{N}[
 \texttt b(X_{T}^{i}) -\texttt b(x_0) -
 (\mathbb E(\texttt b(X_{T}^{i})) -\texttt b(x_0))]\right| >\mathfrak u\right]
\end{eqnarray*}
with
\begin{displaymath}
\mathfrak u =
\log\left(\frac{\mathfrak c\mathfrak e}{
\overline\alpha_HT^{2H}\|\varphi\|_{\infty}\|\psi\|_{\infty}}\right)
\frac{\mathfrak e}{\|\psi\|_{\infty}} -
|\mathbb E(\texttt b(X_T)) -\texttt b(x_0)| > 0
\quad {\rm by}\quad (\ref{deviation_bound_1}).
\end{displaymath}
By the Bienaym\'e-Tchebychev inequality, and since $X^1,\dots,X^N$ are i.i.d. processes,
\begin{displaymath}
\mathbb P(\Delta_{N}^{\mathfrak e})
\leqslant
\frac{1}{\mathfrak u^2N^2}{\rm var}\left(\sum_{i = 1}^{N}
(\texttt b(X_{T}^{i}) -\texttt b(x_0))\right)
\leqslant
\frac{1}{\mathfrak u^2N}
\mathbb E(|\texttt b(X_T) -\texttt b(x_0)|^2).
\end{displaymath}
Therefore,
\begin{displaymath}
\mathbb P(\Delta_{N}^{c})
\leqslant
\frac{\mathfrak c_{\ref{deviation_bound}}}{N}
\quad {\rm with}\quad
\mathfrak c_{\ref{deviation_bound}} =
\frac{\|b^2\|_{f}^{2}}{\mathfrak e^2} +
\frac{1}{\mathfrak u^2}
\mathbb E(|\texttt b(X_T) -\texttt b(x_0)|^2)
\hspace{0.75cm}\Box
\end{displaymath}
%


%
\subsection{Proof of Proposition \ref{ACI_approximation_LS_Young}}\label{section_proof_ACI_approximation_LS_Young}
First, consider
\begin{displaymath}
U_N =\frac{1}{N}\sum_{i = 1}^{N}Z^i,
\end{displaymath}
where $Z^1,\dots,Z^N$ are defined by
\begin{displaymath}
Z^i =\frac{1}{T}\int_{0}^{T}b(X_{s}^{i})\sigma(X_{s}^{i})\delta B_{s}^{i},
\quad\forall i\in\{1,\dots,N\}.
\end{displaymath}
On the one hand, by the (usual) law of large numbers,
\begin{displaymath}
D_N =
\frac{1}{NT}\sum_{i = 1}^{N}\int_{0}^{T}b(X_{s}^{i})^2ds
\xrightarrow[N\rightarrow\infty]{\mathbb P}
\mathbb E\left(\frac{1}{T}\int_{0}^{T}b(X_s)^2ds\right) =\|b\|_{f}^{2} > 0,
\end{displaymath}
By the (usual) central limit theorem,
\begin{displaymath}
\sqrt NU_N\xrightarrow[N\rightarrow\infty]{\mathcal D}\mathcal N(0,{\rm var}(Z))
\quad {\rm with}\quad
Z =\frac{1}{T}\int_{0}^{T}b(X_s)\sigma(X_s)\delta B_s.
\end{displaymath}
Then, since $dX_{t}^{i} =\theta_0b(X_{t}^{i})dt +\sigma(X_{t}^{i})dB_{t}^{i}$ for every $i\in\{1,\dots,N\}$, and by Slutsky's lemma,
\begin{displaymath}
\sqrt N(\widehat\theta_N -\theta_0) =
\sqrt N\frac{U_N}{D_N}
\xrightarrow[N\rightarrow\infty]{\mathcal D}
\mathcal N\left(0,\frac{{\rm var}(Z)}{\|b\|_{f}^{4}}\right).
\end{displaymath}
On the other hand, let us introduce a random variable $Y_{N}^{*}$ converging towards ${\rm var}(Z) = \mathbb E(Z^2)$. Precisely,
\begin{displaymath}
Y_{N}^{*} =
\frac{1}{NT^2}
\sum_{i = 1}^{N}\left(
\alpha_H\int_{0}^{T}\int_{0}^{T}
\pi(X_{s}^{i})\pi(X_{t}^{i})|t - s|^{2H - 2}dsdt + R_i\right)
\end{displaymath}
where, for every $i\in\{1,\dots,N\}$,
\begin{eqnarray*}
 R_i & = &
 \alpha_{H}^{2}\int_{[0,T]^2}\int_{0}^{v}\int_{0}^{u}
 |u -\overline u|^{2H - 2}|v -\overline v|^{2H - 2}\\
 & &
 \hspace{3cm}\times
 \varphi(X_{v}^{i})\varphi(X_{u}^{i})\exp\left(\theta_0\left(\int_{\overline u}^{v}\psi(X_{s}^{i})ds +
 \int_{\overline v}^{u}\psi(X_{s}^{i})ds\right)\right)
 d\overline ud\overline vdudv.
\end{eqnarray*}
By the law of large numbers, and by \cite{BHOZ08}, Theorem 3.11.1,
\begin{displaymath}
Y_{N}^{*}
\xrightarrow[N\rightarrow\infty]{\mathbb P}
{\rm var}(Z) > 0.
\end{displaymath}
So, by Slutsky's lemma,
\begin{displaymath}
\sqrt{
\frac{ND_{N}^{2}}{Y_{N}^{*}}}
(\widehat\theta_N -\theta_0)
\xrightarrow[N\rightarrow\infty]{\mathcal D}\mathcal N(0,1)
\end{displaymath}
and then, for every $x\in\mathbb R_+$,
\begin{displaymath}
\mathbb P\left(\sqrt{
\frac{ND_{N}^{2}}{Y_{N}^{*}}}\cdot
|\widehat\theta_N -\theta_0| > x\right)
\xrightarrow[N\rightarrow\infty]{} 2(1 -\phi(x)).
\end{displaymath}
Now, for $N$ large enough, consider
\begin{displaymath}
c_N =\sqrt{\frac{ND_{N}^{2}}{Y_N}}
\quad {\rm and}\quad
c_{N}^{*} =\sqrt{\frac{ND_{N}^{2}}{Y_{N}^{*}}}.
\end{displaymath}
Since $\psi\leqslant 0$ and $\theta_0 > 0$, $Y_{N}^{*}\leqslant Y_N$, and then $c_{N}^{*}\geqslant c_N$. Moreover,
\begin{eqnarray*}
 |\overline\theta_N -\widehat\theta_N|\mathbf 1_{\Delta_N} & = &
 |\Theta_N(R_N) -\Theta_N(\theta_0 - I_N)|\mathbf 1_{\Delta_N}\\
 & \leqslant &
 \mathfrak c|R_N - (\theta_0 - I_N)|\mathbf 1_{\Delta_N}
 \leqslant
 \mathfrak c|\overline\theta_N -\widehat\theta_N|\mathbf 1_{\Delta_N} +
 \mathfrak c|\widehat\theta_N -\theta_0|\mathbf 1_{\Delta_N},
\end{eqnarray*}
and since $\mathfrak c\in (0,1)$,
\begin{equation}\label{ACI_approximation_LS_Young_1}
|\overline\theta_N -\widehat\theta_N|\mathbf 1_{\Delta_N}
\leqslant
\frac{\mathfrak c}{1 -\mathfrak c}
|\widehat\theta_N -\theta_0|\mathbf 1_{\Delta_N}.
\end{equation}
So, since
\begin{displaymath}
\mathfrak c < 1 -\varepsilon
\quad {\rm with}\quad
\varepsilon =
\frac{u_{1 -\frac{\lambda\alpha}{2}}}{
u_{1 -\frac{\alpha}{2}}}\in (0,1),
\end{displaymath}
for any $x\in\mathbb R_+$,
\begin{eqnarray*}
 \mathbb P(c_N|\overline\theta_{N}^{\mathfrak c} -\theta_0| > x)
 & \leqslant &
 \mathbb P(c_{N}^{*}|\overline\theta_{N}^{\mathfrak c} -\theta_0| > x)\\
 & \leqslant &
 \mathbb P(c_{N}^{*}|\overline\theta_{N}^{\mathfrak c} -\widehat\theta_N| > (1 -\varepsilon)x) +
 \mathbb P(c_{N}^{*}|\widehat\theta_N -\theta_0| >\varepsilon x)\\
 & \leqslant &
 \mathbb P(\Delta_{N}^{c}) +
 \mathbb P(\{c_{N}^{*}|\widehat\theta_N -\theta_0| > (1/\mathfrak c - 1)(1 -\varepsilon)x\}\cap\Delta_N) +
 \mathbb P(c_{N}^{*}|\widehat\theta_N -\theta_0| >\varepsilon x)\\
 & \leqslant &
 \mathbb P(\Delta_{N}^{c}) + 2\mathbb P(c_{N}^{*}|\widehat\theta_N -\theta_0| >\varepsilon x).
\end{eqnarray*}
Therefore,
\begin{displaymath}
\mathbb P(c_N|\overline\theta_{N}^{\mathfrak c} -\theta_0|\leqslant x)
\geqslant 1 -\mathbb P(\Delta_{N}^{c}) -
2\mathbb P(c_{N}^{*}|\widehat\theta_N -\theta_0| >\varepsilon x).
\end{displaymath}
Since, by Proposition \ref{deviation_bound},
\begin{displaymath}
\lim_{N\rightarrow\infty}
\mathbb P(\Delta_{N}^{c}) = 0
\quad\textrm{under the condition (\ref{deviation_bound_1})};
\end{displaymath}
for every $x\in\mathbb R_+$,
\begin{eqnarray*}
 \lim_{N\rightarrow\infty}
 \mathbb P(c_N|\overline\theta_{N}^{\mathfrak c} -\theta_0|\leqslant x)
 & \geqslant &
 1 - 2\underbrace{
 \lim_{N\rightarrow\infty}\mathbb P(c_{N}^{*}|\widehat\theta_N -\theta_0| >\varepsilon x)}_{
 = 2(1 -\phi(\varepsilon x))} = 4\phi(\varepsilon x) - 3.
\end{eqnarray*}
Then, since $4\phi(\varepsilon x) - 3 = 1 - 2\lambda\alpha$ if and only if $x =\varepsilon^{-1}u_{1 -\frac{\lambda\alpha}{2}}$, and by the definition of $\varepsilon$,
\begin{eqnarray*}
 \lim_{N\rightarrow\infty}
 \mathbb P(c_N|\overline\theta_{N}^{\mathfrak c} -
 \theta_0|\leqslant u_{1 -\frac{\alpha}{2}})
 & = &
 \lim_{N\rightarrow\infty}
 \mathbb P(c_N|\overline\theta_{N}^{\mathfrak c} -
 \theta_0|\leqslant\varepsilon^{-1}u_{1 -\frac{\lambda\alpha}{2}})\\
 & \geqslant &
 1 - 2\lambda\alpha
 \hspace{0.75cm}\Box
\end{eqnarray*}
%


%
\subsection{Proof of Proposition \ref{risk_bound_approximation_LS_Young}}
First, since $dX_{t}^{i} =\theta_0b(X_{t}^{i})dt +\sigma(X_{t}^{i})dB_{t}^{i}$ for every $i\in\{1,\dots,N\}$,
\begin{displaymath}
\widehat\theta_N =
\theta_0 +\frac{U_N}{D_N}
\quad {\rm with}\quad
U_N =
\frac{1}{NT}\sum_{i = 1}^{N}\int_{0}^{T}b(X_{s}^{i})\sigma(X_{s}^{i})\delta B_{s}^{i}.
\end{displaymath}
Let us establish a non-asymptotic risk bound on the auxiliary non-computable estimator $\widehat\theta_{N}^{\mathfrak d} =\widehat\theta_N\mathbf 1_{D_N >\mathfrak d}$ of $\theta_0$. Observe that
\begin{eqnarray*}
 \mathbb E(U_{N}^{2}) & = &
 \frac{1}{N^2T^2}
 \sum_{i = 1}^{N}{\rm var}
 \left(\int_{0}^{T}b(X_{s}^{i})\sigma(X_{s}^{i})\delta B_{s}^{i}\right)\\
 & = &
 \frac{\mathfrak c_1}{N}
 \quad {\rm with}\quad
 \mathfrak c_1 =\frac{1}{T^2}
 {\rm var}\left(\int_{0}^{T}b(X_s)\sigma(X_s)\delta B_s\right).
\end{eqnarray*}
Moreover, a computation along the same lines as in (\ref{deviation_bound_2}) guarantees that
\begin{displaymath}
\mathbb P(D_N <\mathfrak d)
\leqslant
\frac{\|b^2\|_{f}^{2}}{\mathfrak d^2N}.
\end{displaymath}
So,
\begin{eqnarray}
 \mathbb E(|\widehat\theta_{N}^{\mathfrak d} -\theta_0|^2)
 & \leqslant &
 \mathbb E(|\widehat\theta_N -\theta_0|^2\mathbf 1_{D_N\geqslant\mathfrak d}) +
 \theta_{0}^{2}\mathbb P(D_N <\mathfrak d)
 \nonumber\\
 \label{risk_bound_approximation_LS_Young_1}
 & \leqslant &
 \frac{1}{\mathfrak d^2}\left(\mathbb E(U_{N}^{2}) +
 \theta_{0}^{2}\frac{\|b^2\|_{f}^{2}}{N}\right)
 \leqslant\frac{\mathfrak c_2}{N}
 \quad {\rm with}\quad
 \mathfrak c_2 =\frac{1}{\mathfrak d^2}(\mathfrak c_1 +\theta_{0}^{2}\|b^2\|_{f}^{2}).
\end{eqnarray}
Now, Inequality (\ref{ACI_approximation_LS_Young_1}) leads to
\begin{eqnarray*}
 |\overline\theta_N -\theta_0|\mathbf 1_{\Delta_N}
 & \leqslant &
 |\overline\theta_N -\widehat\theta_N|\mathbf 1_{\Delta_N} +
 |\widehat\theta_N -\theta_0|\mathbf 1_{\Delta_N}\\
 & \leqslant &
 \mathfrak c_3
 |\widehat\theta_N -\theta_0|\mathbf 1_{\Delta_N}
 \quad {\rm with}\quad
 \mathfrak c_3 =
 \frac{1}{1 -\mathfrak c} > 1.
\end{eqnarray*}
So,
\begin{eqnarray*}
 |\overline\theta_{N}^{\mathfrak c,\mathfrak d} -\theta_0|
 & = &
 |\overline\theta_N -\theta_0|\mathbf 1_{\{D_N\geqslant\mathfrak d\}\cap\Delta_N} +
 |\theta_0|\mathbf 1_{(\{D_N\geqslant\mathfrak d\}\cap\Delta_N)^c}\\
 & \leqslant &
 \mathfrak c_3
 |\widehat\theta_N -\theta_0|\mathbf 1_{\{D_N\geqslant\mathfrak d\}\cap\Delta_N} +
 |\theta_0|(\mathbf 1_{D_N <\mathfrak d} +\mathbf 1_{\Delta_{N}^{c}})\\
 & \leqslant &
 \mathfrak c_3
 (|\widehat\theta_N -\theta_0|\mathbf 1_{D_N\geqslant\mathfrak d} +
 |\theta_0|\mathbf 1_{D_N <\mathfrak d}) +
 |\theta_0|\mathbf 1_{\Delta_{N}^{c}}\\
 & = &
 \mathfrak c_3
 |\widehat\theta_{N}^{\mathfrak d} -\theta_0| +
 |\theta_0|\mathbf 1_{\Delta_{N}^{c}}.
\end{eqnarray*}
Therefore, by Inequality (\ref{risk_bound_approximation_LS_Young_1}) together with Proposition \ref{deviation_bound},
\begin{displaymath}
\mathbb E(|\overline\theta_{N}^{\mathfrak c,\mathfrak d} -\theta_0|^2)\leqslant
\frac{2\mathfrak c_1\mathfrak c_{3}^{2}}{N} +
2\theta_{0}^{2}\mathbb P(\Delta_{N}^{c})
\lesssim\frac{1}{N}
\hspace{0.75cm}\Box
\end{displaymath}
%


%
\subsection{Proof of Proposition \ref{existence_uniqueness_R_N_rough}}\label{section_proof_existence_uniqueness_R_N_rough}
First, since $\alpha_H = H(2H - 1) < 0$ and $\psi\leqslant 0$, for every $i\in\{1,\dots,N\}$, $t\in [0,T]$ and $\theta\in\mathbb R_+$,
\begin{equation}\label{existence_uniqueness_R_N_rough_2}
\Lambda_{t}^{i}(\theta) = -Ht^{2H - 1} +
\alpha_H\int_{0}^{t}\left(1 -\exp\left(\theta\int_{s}^{t}\psi(X_{u}^{i})du\right)\right)|t - s|^{2H - 2}ds
\leqslant 0.
\end{equation}
Then, since $\varphi\leqslant 0$ and $I_N\geqslant 0$,
\begin{displaymath}
\widetilde\Theta_N(r) =\frac{1}{NTD_N}\sum_{i = 1}^{N}
\int_{0}^{T}\varphi(X_{t}^{i})\Lambda_{t}^{i}(r + I_N)dt\geqslant 0,
\quad\forall r\in\mathbb R_+.
\end{displaymath}
Now, for any $r,\overline r\in\mathbb R_+$,
\begin{eqnarray*}
 \widetilde\Theta_N(\overline r) -\widetilde\Theta_N(r) & = &
 -\frac{\alpha_H}{NTD_N}\sum_{i = 1}^{N}
 \int_{0}^{T}\int_{0}^{t}|t - s|^{2H - 2}\varphi(X_{t}^{i})
 \exp\left(I_N\int_{s}^{t}\psi(X_{u}^{i})du\right)\\
 & &
 \hspace{4cm}\times
 \left(\exp\left(\overline r\int_{s}^{t}\psi(X_{u}^{i})du\right) -
 \exp\left(r\int_{s}^{t}\psi(X_{u}^{i})du\right)\right)dsdt.
\end{eqnarray*}
So, exactly as in the proof of Proposition \ref{existence_uniqueness_R_N_Young},
\begin{eqnarray*}
 |\widetilde\Theta_N(\overline r) -\widetilde\Theta_N(r)|
 & \leqslant &
 \frac{|\alpha_H|\cdot\|\varphi\|_{\infty}}{NTD_N}
 e^{\|\psi\|_{\infty}|I_N|T}
 \sum_{i = 1}^{N}\int_{0}^{T}\int_{0}^{t}|t - s|^{2H - 2}\sup_{x\in\mathbb R_-}e^x\left|
 (\overline r - r)\int_{s}^{t}\psi(X_{u}^{i})du\right|dsdt\\
 & \leqslant &
 \overline\alpha_H\|\varphi\|_{\infty}\|\psi\|_{\infty}
 T^{2H}\frac{M_N}{D_N}|\overline r - r|
 \leqslant\mathfrak c|\overline r - r|.
\end{eqnarray*}
Therefore, $\widetilde\Theta_N$ is a contraction from $\mathbb R_+$ into itself, and then $R_N$ exists and is unique by Picard's fixed-point theorem.$\hspace{0.75cm}\Box$
%


%
\subsection{Proof of Proposition \ref{ACI_approximation_LS_rough}}
First, as in the proof of Proposition \ref{ACI_approximation_LS_Young},
\begin{displaymath}
\sqrt N(\widehat\theta_N -\theta_0) =
\sqrt N\frac{U_N}{D_N}
\xrightarrow[N\rightarrow\infty]{\mathcal D}
\mathcal N\left(0,\frac{{\rm var}(Z)}{\|b\|_{f}^{4}}\right)
\quad {\rm with}\quad
Z =\frac{1}{T}\int_{0}^{T}\pi(X_{s}^{i})\delta B_{s}^{i}.
\end{displaymath}
Consider
\begin{eqnarray*}
 \mathfrak Y_{N}^{*} & = &
 \frac{1}{N}\sum_{i = 1}^{N}\left(
 \frac{1}{T}\int_{0}^{T}
 \pi(X_{s}^{i})\delta B_{s}^{i}\right)^2\\
 & = &
 \frac{1}{NT^2}\sum_{i = 1}^{N}\left(
 \int_{0}^{T}\pi(X_{s}^{i})dB_{s}^{i} +\int_{0}^{T}\varphi(X_{t}^{i})\Lambda_{t}^{i}(\theta_0)dt\right)^2
 \quad {\rm by}\quad (\ref{relationship_Skorokhod_rough_1}).
\end{eqnarray*}
By the law of large numbers,
\begin{displaymath}
\mathfrak Y_{N}^{*}
\xrightarrow[N\rightarrow\infty]{\mathbb P}
{\rm var}(Z) > 0.
\end{displaymath}
So, by Slutsky's lemma,
\begin{displaymath}
\sqrt{
\frac{ND_{N}^{2}}{\mathfrak Y_{N}^{*}}}
(\widehat\theta_N -\theta_0)
\xrightarrow[N\rightarrow\infty]{\mathcal D}\mathcal N(0,1)
\end{displaymath}
and then, for every $x\in\mathbb R_+$,
\begin{displaymath}
\mathbb P\left(\sqrt{
\frac{ND_{N}^{2}}{\mathfrak Y_{N}^{*}}}\cdot
|\widehat\theta_N -\theta_0| > x\right)
\xrightarrow[N\rightarrow\infty]{} 2(1 -\phi(x)).
\end{displaymath}
Now, recall that $\mathfrak Y_N$ was introduced in the statement of this proposition, and - for $N$ large enough - consider
\begin{displaymath}
c_N =\sqrt{\frac{ND_{N}^{2}}{\mathfrak Y_N}}
\quad {\rm and}\quad
c_{N}^{*} =\sqrt{\frac{ND_{N}^{2}}{\mathfrak Y_{N}^{*}}}.
\end{displaymath}
For any $i\in\{1,\dots,N\}$, since $\theta_0\leqslant\theta_{\max}$,
\begin{eqnarray*}
 \left|\int_{0}^{T}\pi(X_{s}^{i})dB_{s}^{i}\right|
 & = &
 \left|\int_{0}^{T}\frac{\pi(X_{s}^{i})}{\sigma(X_{s}^{i})}
 (\theta_0b(X_{s}^{i})ds +\sigma(X_{s}^{i})dB_{s}^{i} -\theta_0b(X_{s}^{i})ds)\right|\\
 & \leqslant &
 \left|\int_{0}^{T}b(X_{s}^{i})dX_{s}^{i}\right| +
 \theta_{\max}\int_{0}^{T}b(X_{s}^{i})^2ds.
\end{eqnarray*}
Moreover, by Inequality (\ref{existence_uniqueness_R_N_rough_2}), and since $\alpha_H < 0$ and $\varphi,\psi\leqslant 0$,
\begin{displaymath}
0\leqslant
\int_{0}^{T}\varphi(X_{t}^{i})\Lambda_{t}^{i}(\theta_0)dt\leqslant
\int_{0}^{T}\varphi(X_{t}^{i})\Lambda_{t}^{i}(\theta_{\max})dt.
\end{displaymath}
So,
\begin{displaymath}
\mathfrak Y_{N}^{*}\leqslant\mathfrak Y_N,
\quad\textrm{leading to}\quad
c_{N}^{*}\geqslant c_N.
\end{displaymath}
Therefore, as in the proof of Proposition \ref{ACI_approximation_LS_Young},
\begin{displaymath}
\mathbb P(c_N|\overline\theta_{N}^{\mathfrak c} -\theta_0|\leqslant x)
\geqslant 1 -\mathbb P(\Delta_{N}^{c}) -
2\mathbb P(c_{N}^{*}|\widehat\theta_N -\theta_0| >\varepsilon x)
\quad {\rm with}\quad
\varepsilon =
\frac{u_{1 -\frac{\lambda\alpha}{2}}}{
u_{1 -\frac{\alpha}{2}}}.
\end{displaymath}
In conclusion, by Proposition \ref{deviation_bound},
\begin{displaymath}
\lim_{N\rightarrow\infty}
\mathbb P(c_N|\overline\theta_{N}^{\mathfrak c} -
\theta_0|\leqslant u_{1 -\frac{\alpha}{2}})
\geqslant 1 - 2\lambda\alpha
\hspace{0.75cm}\Box
\end{displaymath}
%


%

%
\appendix
%


%
\section{Summary of key notations}\label{key_notations_section}
To help the reader in navigating the technical framework of this paper, we collect and summarize the core notation and functions used throughout the proofs. Due to the high number of interacting symbols, this section serves as a centralized reference guide.
\\
\\
{\bf Auxiliary functions.} To express the Malliavin derivatives independently of $B$, we define the following auxiliary functions based on the drift and volatility functions:
\begin{displaymath}
\pi = b\sigma,\quad
\varphi =\sigma\pi' =\sigma(\sigma b' +\sigma'b)
\quad {\rm and}\quad
\psi =\frac{1}{\sigma}(\sigma b' -\sigma'b).
\end{displaymath}
{\bf Empirical quantities.} Several empirical quantities are involved in the definition of our fixed-point maps:
\begin{displaymath}
D_N =\frac{1}{NT}
\sum_{i = 1}^{N}\int_{0}^{T}b(X_{s}^{i})^2ds,\quad
I_N =\frac{1}{NTD_N}
\sum_{i = 1}^{N}(\mathtt{b}(X_{T}^{i}) - \mathtt{b}(x_0))
\quad (\texttt b' = b)\quad {\rm and}\quad
M_N = e^{\|\psi\|_{\infty}|I_N|T}.
\end{displaymath}
{\bf Fixed-point maps.} Depending on the Hurst parameter $H$, we define the random maps whose fixed-points are involved in the definition of our main estimators. Both maps heavily rely on the quantities $D_N$, $I_N$ and $M_N$. For the regime $H > 1/2$, the computable estimator of $\theta_0$ depends on the fixed-point of the map
\begin{displaymath}
\Theta_N(\cdot) :=
-\frac{\alpha_H}{NTD_N}
\sum_{i = 1}^{N}\int_{0}^{T}\int_{0}^{t}
\varphi(X_{t}^{i})\exp\left(
(\cdot + I_N)\int_{s}^{t}\psi(X_{u}^{i})du\right)|t - s|^{2H - 2}dsdt.
\end{displaymath}
For the regime $H\in (1/3,1/2]$, the computable estimator of $\theta_0$ depends on the fixed-point of the map
\begin{displaymath}
\widetilde\Theta_N(\cdot) :=
\frac{1}{NTD_N}
\sum_{i = 1}^{N}\int_{0}^{T}\varphi(X_{t}^{i})
\left[-Ht^{2H - 1} +\alpha_H\int_{0}^{t}\left(
1 -\exp\left((\cdot + I_N)\int_{s}^{t}\psi(X_{u}^{i})du\right)\right)
|t - s|^{2H - 2}ds\right]dt.
\end{displaymath}
{\bf The "stability" event.} Finally, to ensure the existence and uniqueness of the fixed-point in both regimes, the same underlying condition appears. This naturally leads to the event
\begin{displaymath}
\Delta_N :=
\left\{T^{2H}\frac{M_N}{D_N}\leqslant\frac{\mathfrak c}{
\overline\alpha_H\|\varphi\|_{\infty}\|\psi\|_{\infty}}\right\}
\quad {\rm with}\quad
\overline\alpha_H =\frac{|\alpha_H|}{2H(2H + 1)}.
\end{displaymath}
\end{document}